\documentclass[12pt]{article}
\usepackage{graphicx} 
\usepackage[a4paper, margin = 2.33cm, centering]{geometry}
\usepackage{enumerate}
\usepackage{subfigure}
\usepackage{algorithm}
\usepackage{algpseudocode}
\usepackage{authblk}
\usepackage{natbib}
\usepackage{hyperref}       
\usepackage{appendix}
\usepackage{xcolor}

\usepackage{setspace}

\usepackage{lscape}

\usepackage{lineno}
\usepackage{arydshln} 
\usepackage{multirow}
\usepackage{latexsym,mathtools,mathrsfs,amssymb,amsmath,amsbsy,amsopn,amstext,amsthm,amsxtra,bm,bbm,calc,ifpdf,cancel,extarrows,manfnt,cancel}

\newcommand{\1}{\mathbbm{1}}				
\renewcommand{\P}{\mathbb{P}}				
\newcommand{\E}{\mathbb{E}}					
\newcommand{\R}{\mathbf{R}}	
 
\newcommand{\argmin}{\operatornamewithlimits{arg\,min}}

\newcommand{\Q}{\mathbf{Q}}
\newcommand{\A}{\mathbf{A}}
\newcommand{\K}{\mathcal{K}}
\newcommand{\T}{\bm{\Theta}}
\newcommand{\J}{\mathcal{J}}
\newcommand{\N}{\mathcal{N}}

\newcommand{\pto}{\overset{\mathbb P}{\longrightarrow}} 
\newcommand{\abs}[1]{\left\lvert#1\right\rvert}		

\newcommand{\hta}{\widehat{\theta}_{j,a}}
\newcommand{\bta}{\overline{\theta}_{j,a}}
\newcommand{\btone}{\overline{\theta}_{j,1}^{(\widehat{\mathbf{A}})}}
\newcommand{\bttwo}{\overline{\theta}_{j,2}^{(\widehat{\mathbf{A}})}}
\newcommand{\tone}{\theta^0_{j,1}}
\newcommand{\ttwo}{\theta^0_{j,2}}
\newcommand{\malpha}{\boldsymbol{\alpha}}
\newcommand{\vtheta}{\bm{\theta}}
\newcommand{\bmq}{\boldsymbol{q}}

\newtheorem{Lem}{Lemma}
\newtheorem{Thm}{Theorem}

\newtheorem{Eg}{Example}
\newtheorem{Ass}{Assumption}

\newtheorem{Prop}{Proposition}

\title{{Consistency Theory} of General Nonparametric Classification Methods in Cognitive Diagnosis}

\author{Chengyu Cui$^1$, Yanlong Liu$^2$, and  Gongjun Xu$^1$\footnote{Cui and Liu are co-first authors}}
\affil{University of Michigan$^1$ and University of Chicago$^2$} 
\author{}
\date{}

\begin{document}
\maketitle

\begin{abstract}
    Cognitive diagnosis models have been popularly used in fields such as education, psychology, and social sciences. While parametric likelihood estimation is a prevailing method for fitting cognitive diagnosis models, nonparametric methodologies are attracting increasing attention due to their ease of implementation and robustness, particularly when sample sizes are relatively small. However, existing {consistency} results of the nonparametric estimation methods often rely on certain restrictive conditions, which may not be easily satisfied in practice. 
    In this article, the {consistency theory for} the general nonparametric classification method is reestablished under weaker and more practical conditions.    

    \vspace{2pt}
\end{abstract}

\section{Introduction}\label{sec:intro}

Cognitive diagnosis models (CDMs), also known as diagnostic classification models (DCMs), are a popular family of discrete latent variable models employed in diagnostic assessments to provide detailed information about subjects' latent attributes based on their responses to designed diagnostic items. For instance, in educational testing, these latent attributes might indicate if a subject has mastered certain skills or not ~\citep{delatorre2011, henson2009defining, junker2001}; in psychiatric diagnosis, the latent attributes might signal the presence or absence of certain mental disorders~\citep{templin2006, delaTorre2018analysis}. 

Parametric models for cognitive diagnosis have been developed and widely applied in practice. Popular examples include the deterministic input, noisy ``and'' gate (DINA) model~\citep{junker2001}, the deterministic input, noisy ``or'' gate (DINO) model~\citep{templin2006}, the general diagnostic model (GDM; ~\citeauthor{davier2008}, 2008), the reduced reparameterized unified model (Reduced RUM;~\citeauthor{hartz2002bayesian}, 2002), the log-linear CDM (LCDM;~\citeauthor{henson2009defining}, 2009), and the generalized DINA model (GDINA;~\citeauthor{delatorre2011}, 2011). In conventional settings with a fixed number of items ($J$) and a large number of subjects ($N$), the latent attributes are often viewed as random variables. 
The corresponding CDMs can thus be viewed as a family of finite mixture models, where each subject's latent attribute profile $\malpha_i$ behaves as a discrete random variable following a categorical distribution. From this perspective, the estimation often takes place through the maximization of the marginal likelihood, relying on methods such as the expectation-maximization algorithm 
\citep{delatorre2011, davier2008, dibello2007review}. However, the maximum likelihood-based approach often necessitates sufficiently large assessments to guarantee the reliability of the item parameter estimation, and it may either produce inaccurate estimates with small sample sizes or
suffer from high computational costs~\citep{chiu2019consistency, chiu2018cognitive}.
Moreover, the parametric CDMs involve certain parametric assumptions about the item response functions, which may raise concerns about the validity of the assumed model and the underlying process~\citep{chiu2013nonparametric}.

As an alternative, researchers have explored nonparametric cognitive diagnosis methods~\citep{chiu2009cluster,chiu2019nonparametric}.
Instead of modeling the item response functions parametrically, the nonparametric methods aim to directly categorize subjects into latent groups by minimizing certain distance measure between a subject's observed item responses and some expected ``centers" of the latent groups. Two popular examples of nonparametric cognitive diagnosis methods 
include the nonparametric classification (NPC) method~\citep{chiu2013nonparametric} and its generalization, the general NPC (GNPC) method~\citep{chiu2018cognitive}. The GNPC method, in particular, has received increasing attention in recent years due to its effectiveness in handling complex CDMs and its good performance for sample sizes~\citep{chiu2021advances,wang2023general, chandia2023nonparametric,Ma2023}. The algorithms of the NPC and GNPC methods are straightforward to implement and require minimal computational resources, making them highly appealing for practical applications.

Theoretical properties of the nonparametric methods have also been explored in the literature. 
Under some regularity conditions, the NPC estimators of the subjects' latent attribute profiles have been shown to be statistically consistent for certain  CDMs, including  DINA and Reduced RUM ~\citep{wang2015}, and a similar {consistency theory} for the GNPC estimator has also been established~\citep{chiu2019consistency}. 
However, the current theoretical guarantees for these nonparametric methods depend on relatively stringent assumptions. In the case of the NPC method, the assumptions associated with the ideal binary responses might oversimplify the underlying diagnostic process and thus be challenging to fulfill when dealing with complex underlying cognitive diagnosis models, such as the GDINA model and other general CDMs~\citep{chiu2018cognitive}. 
{Although the GNPC method addresses the oversimplification issue of the NPC method, its consistency depends on a key assumption that consistent initial estimators of the latent attribute profiles are available. For instance, Theorem 1 in \cite{chiu2019consistency} provides theoretical guarantees for the GNPC estimators, assuming an initialization that consistently estimates the ground truth latent memberships. Similarly, Theorems 1--3 in \cite{Ma2023} require consistent estimation of latent memberships from a calibration dataset to establish their consistency results.}
The assumption that consistent initial estimators of latent attribute profiles can be obtained or that a calibration dataset is available may be overly restrictive in practice, and the consistency of the GNPC method in more realistic settings remains an open problem.

In this article, we establish the {consistency} for the GNPC method using different theoretical techniques, without requiring the previous assumption on initial consistent estimators or calibration datasets. Our analysis covers both the original GNPC method in~\cite{chiu2019consistency}  and a modified version of the GNPC method in~\cite{Ma2023}.
{
We establish {finite-sample error bounds for latent attributes of general nonparametric methods} as well as the uniform consistency of the item parameters. 
We would like to clarify that the main contribution of this work lies in the theoretical analysis of the GNPC and modified GNPC methods. For the implementation of these methods, we recommend utilizing the algorithms proposed in the literature \citep{chiu2018cognitive,chiu2019consistency,Ma2023}, which have demonstrated the effectiveness of GNPC methods via extensive simulation studies and real data examples.
} 


The rest of the paper is organized as follows.
Section~\ref{sec:setup} provides a brief review of cognitive diagnostic models and discusses the limitations in the existing consistency results. Section~\ref{sec:estimator}  establishes consistency results of the GNPC methods. {In Section~\ref{sec:simu}, we provide a simulation study to illustrate our theoretical results.} Section~\ref{sec:discussion} gives some further discussions, and   the supplementary material provides the proofs for the main results.

\section{Model Setup and Nonparametric Methods}\label{sec:setup}


This work focuses on CDMs for multivariate binary data, which are commonly encountered in educational assessments (correct/wrong answers) and social science survey responses (yes/no responses) \citep{von2019handbook}. For $N$ subjects and $J$ items, the observed data is an $N\times J$ binary matrix $\R=(R_{i,j})$, where $R_{i,j}=1$ or $0$ denotes whether the $i$-th subject gives a positive response to the $j$-th item. Consider $K$ binary latent attributes. Let the row vector \(\malpha_i = (a_{i,1},\dots,a_{i,K})\) represent the latent attribute profile for the \(i\)-th subject, where
 \(a_{i,k} = 1\) or \(0\) indicates the presence or absence, respectively, of the \(k\)-th attribute for the \(i\)-th individual.
We further use an \(N \times K\) binary matrix,
\(\A = (a_{i,k})\in \{0,1\}^{N\times K}\), to represent the latent attribute profiles for all $N$ subjects.  

To capture the dependence relationship between items and the latent attributes of subjects, a design matrix called the Q-matrix~\citep{tatsuoka1985} is employed. The Q-matrix encodes how the $J$ items depend on the $K$ latent attributes. Specifically, $\mathbf{Q}=(q_{j,k})\in\{0,1\}^{J\times K}$, where $q_{j,k}=1$ or $0$ indicates whether the $j$-th test item depends on the $k$-th latent attribute, and we denote the $j$-th item's Q-matrix vector as $\bm{q}_j= (q_{j,1},\dots,q_{j,K})$. 

For an integer $m$, we denote $[m] =\{1,\cdots, m\}$ and for a set ${\cal A}$, we denote  its cardinality by $|{\cal A}|$. We denote $\theta_{j,\malpha}= \P(R_{i,j}=1|\malpha_i = \malpha)$ {for any $i\in [N]$}, $ j\in [J]$ and  $\malpha\in\{0,1\}^K$, and let 
$\T =\{\theta_{j,\malpha};j\in[J],\malpha\in\{0,1\}^K\}$.
We assume each response $R_{i,j}$ follows a Bernoulli distribution with parameter $\theta_{j,\malpha_i}$  and the responses are independent with each other conditional on the latent attribute profiles $\A$ and the structure loading matrix $\Q$. In summary, the data generative process aligns with the following  latent class model: 
\begin{align}
    \P(\R\mid\A,\T) = \prod_{i=1}^N \prod_{j=1}^J\P(R_{i,j}\mid\malpha_{i},\theta_{j,\malpha_i})=\prod_{i=1}^N \prod_{j=1}^J (\theta_{j,\malpha_i})^{R_{i,j}} (1 - \theta_{j,\malpha_i})^{1-R_{i,j}}. \nonumber
\end{align}
To further illustrate the adaptability of the general nonparametric method to the model structures  embedded in CDMs imposed by the structural matrix \(\Q\), we follow the general assumption for the  restricted latent class models ~\citep{chiu2015consistency,xu2017identifiability,Ma2023} that for different attribute profiles $\tilde{\malpha}$ and $\malpha$, 
\begin{align} \label{eq3}\left(\malpha \circ \bm{q}_j = \tilde{\malpha} \circ \bm{q}_j \right) \Longrightarrow \left( \theta_{j,\malpha} = \theta_{j,\tilde{\malpha}}\right),\end{align}
where {$\malpha \circ \bm{q}_j = (a_1  q_{j,1},\dots,a_K  q_{j,K})$ denotes the element-wise product of binary vectors $\malpha=(\alpha_1,\cdots,\alpha_K)$ and $\bm{q}_j$}. This implies that the item response parameter $\theta_{j,\malpha}$ only depends on whether the latent attribute profile  $\malpha$ contains the required attributes $\K_j:=\{k\in[K];q_{j,k}=1 \}$ for item $j$. 
In cognitive diagnostic assessments, the matrix $\Q$ is typically pre-determined by domain experts~\citep{george2015cognitive,junker2001,davier2008}. 
In this work, we assume the Q-matrix $\Q$ is specified, and $(\A,\T)$ are to be estimated from the responses  $\R$.

\paragraph{Parametric CDMs: DINA and DINO models}

For parametric CDMs, the structural matrix $\Q$ imposes various constraints on the item parameters based on different cognitive assumptions. For instance, in the DINA~\citep{junker2001} model, a conjunctive relationship among the attributes is assumed. According to this assumption, for a subject to provide a positive (correct) response to an item, mastery of all the required attributes of the item is necessary. In the DINA model, the ideal response for each item $j\in [J]$ and each latent attribute profile $\malpha=(a_1,\dots,a_K)$ is defined as $$\eta_{j,\malpha}^{\mathrm{DINA}} = \prod_{k=1}^K a_k^{q_{j,k}}.$$

\noindent The DINO~\citep{templin2006} model assumes a disjunctive relationship among attributes, where mastery of at least one of the required attributes for an item is necessary for a subject to be considered capable of providing a positive response. In the DINO model, the ideal response is defined as $$\eta_{j,\malpha}^{\mathrm{DINO}} = 1- \prod_{k=1}^K (1-a_k)^{q_{j,k}}.$$

\noindent The DINA and DINO models further encompass uncertainty by incorporating the slipping and guessing parameters, denoted as \(s_j\) and \(g_j\) for \(j \in [J]\). For each item \(j\), the slipping parameter represents the probability of a capable subject giving a negative response, whereas the guessing parameter signifies the probability of an incapable subject giving a positive response. Specifically, \(s_j = \P(R_{i,j}=0|\eta_{j,\malpha_i}=1)\) and \(g_j = \P(R_{i,j}=1|\eta_{j,\malpha_i}=0)\) for the $i$-th subject. Therefore, in these two restricted latent class models, the parameter $\theta_{j,\malpha}$ can be expressed as $$\theta_{j,\malpha} = (1-s_j)^{\eta_{j,\malpha}} g_j^{1-\eta_{j,\malpha}}.$$


\paragraph{Nonparametric CDMs: NPC and GNPC}
For nonparametric CDMs, the ideal responses described under the DINA and DINO models serve as foundational elements for the nonparametric {classification} analysis. Given a set of 0-1 binary ideal responses, denoted as \(\{\eta_{j,\malpha}\}\), the nonparametric classification (NPC) method, as introduced by~\cite{chiu2013nonparametric}, estimates the subjects' latent attribute profiles as follows. This method utilizes a distance-based algorithm, leveraging observed item responses to categorize subjects into latent groups. The NPC estimator, \(\widehat{\malpha}_i\), for the $i$-th individual's attribute profile, \(\malpha_i\), is expressed as 
\[
\widehat{\malpha}_i = \argmin_{\malpha \in \{0,1\}^K} \sum_{j=1}^J (R_{i,j} - \eta_{j,\malpha})^2.
\]

\noindent In the NPC method, the ideal responses $\eta_{j,\malpha}$ can be based on either the DINA model or the DINO model. However, due to the dependence on these specific model assumptions, which define two extreme relations between $\bmq_j$ and latent attribute profile $\malpha$,  the NPC method may fail to handle complex CDMs, such as the GDINA model, 
and such limitation may lead to misclassifications of the subjects~\citep{chiu2019consistency}. 

To address this issue, the GNPC method~\citep{chiu2018cognitive} offers a solution by considering a more general ideal response that represents a weighted average of the ideal responses from the DINA and DINO models, as in 
\hypertarget{constraint2}{}\begin{align}\eta_{j,\malpha}^{(w)} = w_{j,\malpha} \eta_{j,\malpha}^{\mathrm{DINA}} + (1-w_{j,\malpha}) \eta_{j,\malpha}^{\mathrm{DINO}}. \label{eq4}\end{align}

\noindent The weights are determined by the data; therefore, the proportional influence of $\eta_{j,\malpha}$ and $w_{j,\malpha}$ on the weighted ideal item response is adapted to the complexity of the underlying CDM data generating process. 
The GNPC method can be utilized with any CDM that can be represented as a general CDM, without requiring prior knowledge of the underlying model. To obtain estimates of the weights,~\cite{chiu2018cognitive} proposed minimizing the $L^2$ distance between the responses to item $j$ and the weighted ideal responses $\eta_{j,\malpha}^{(w)}$: \begin{align}d_{j,\malpha} = \sum_{i:{\malpha}_i =\malpha} (R_{i,j} - \eta_{j,\malpha}^{(w)})^2. \nonumber\end{align}
 When $\eta_{j,\malpha}^{\mathrm{DINO}} = \eta_{j,\malpha}^{\mathrm{DINA}}$, this results in $\eta_{j,\malpha}^{(w)} = \eta_{j,\malpha}^{\mathrm{DINA}} = \eta_{j,\malpha}^{\mathrm{DINO}}$, 
 which happens either when $\malpha$ includes all the required latent attributes in $\mathcal{K}_j$, leading to $\eta_{j,\malpha}^{(w)}=1$, or when $\malpha$ does not contain any required attributes, resulting in $\eta_{j,\malpha}^{(w)}=0$. 
 Equivalently, these two  extreme situations can be summarized as the following constraints:
  \begin{equation}    \left({\malpha} \cdot \bm{q}_j = 0 \Longrightarrow \eta_{j,\malpha}^{(w)} = 0\right) 
  \text{ and }  \left({\malpha} \cdot \bm{q}_j = K_j \Longrightarrow \eta_{j,\malpha}^{(w)} = 1\right),\label{eq10-eq}
  \end{equation} 
where {${\malpha} \cdot \bm{q}_j = \sum_{k=1}^K\alpha_{k}q_{j,k}$} denotes the inner product of the two vectors and
\( K_j \) is defined as \( \sum_{k=1}^K q_{j,k} \), representing the number of latent attributes that the $j$-th item depends on. Thus in these two extreme situations, the parameters $\eta_{j,\malpha}^{(w)}$ are known and do not need estimation. 
 In scenarios where $\malpha$ includes only some of the required attributes, 
 $\eta_{j,\malpha}^{(w)}$ need to be estimated, and  in such cases, minimizing $d_{j,\malpha}$ would lead to
 \begin{align}
 \widehat{w}_{j,\malpha} = 1- \overline{R}_{j,\malpha}, \quad \widehat{\eta}_{j,\malpha}^{(w)} = \overline{R}_{j,\malpha},
 \label{eq6}\end{align}
 where $\overline{R}_{j,\malpha} = \sum_{i:{\malpha}_i = \malpha} R_{i,j} / |\{i\in[N];{\malpha}_i = \malpha \}|$, which represents the sample mean of the responses to the $j$-th item for subjects with given latent attribute profile $\malpha$. Since the true latent attribute profiles are unknown, the memberships and the ideal responses will be jointly estimated. Specifically, the optimization problem associated with the GNPC method in \cite{chiu2018cognitive}  aims to minimize the following loss function over the membership $\malpha_i$ and the weights $w_{j,\malpha}$ under the constraints imposed by the given Q-matrix: \begin{align}\sum_{\malpha \in \{0,1\}^K} \sum_{i: {\malpha}_i = \malpha} \sum_{j=1}^J (R_{i,j} - \eta_{j,\malpha}^{(w)})^2, \label{eq7}\end{align}
under constraint~\eqref{eq3}, where $\eta_{j,\malpha}^{(w)}$ is given in \eqref{eq4}.

 {A modified GNPC method was studied by~\cite{Ma2023} under a general framework where the item parameters $\theta_{j,\malpha}$ are treated as a certain ``centroid''. In their framework, the item parameters $\theta_{j,\malpha}$ and latent attributes $\malpha_i$ are obtained by minimizing $L(\A,\T) = \sum_{\malpha\in\{0,1\}^K}\sum_{i:\malpha_i = \malpha} l(\boldsymbol{R}_i,\boldsymbol{\theta}_{\malpha})$ where $l(\boldsymbol{R}_i,\boldsymbol{\theta}_{\malpha})$ is a loss function that measures the distance between the $i$th subject's response vector, $\boldsymbol{R}_i=(R_{i,j}, j=1,\ldots, J)$, and the item parameter vector $\boldsymbol{\theta}_{\malpha} = (\theta_{j,\malpha}, j=1,\ldots, J)$, given a membership $\malpha$. Under their framework, GNPC method can be derived by taking $l(\boldsymbol{R}_i,\boldsymbol{\theta}_{\malpha}) = \sum_{j=1}^J(R_{i,j} - \theta_{j,\malpha})^2$, which leads to minimizing the following loss function}
 \begin{align}\sum_{\malpha \in \{0,1\}^K} \sum_{i: {\malpha}_i = \malpha} \sum_{j=1}^J (R_{i,j} - \theta_{j,\malpha})^2, \label{eq_mGNPC}\end{align}
with respect to $\theta_{j,\malpha}$ and $\malpha$  under constraint \eqref{eq3}. To ensure identifiability, we impose the natural constraint $\theta_{j,\malpha}\ge \theta_{j,\tilde\malpha}$ if $\malpha\succeq\tilde\malpha$. Here $\malpha\succeq\tilde\malpha$ if $\alpha_{k}\ge \tilde\alpha_k$ for all $k\in[K]$.
 
{Note that given the membership $\malpha$, the item parameter $\theta_{j,\malpha}$ that minimizes the loss function~\eqref{eq_mGNPC} takes exactly the form of $\overline{R}_{j,\malpha}$ in \eqref{eq6} for all items and  $\malpha$'s. Inspired by this, as shown in~\cite{Ma2023}, we can see that the solution $(\hat\malpha_i,\hat\eta_{j,\malpha})$ to the original GNPC estimation method in \eqref{eq7} is the same as the solution  $(\hat\malpha_i,\hat\theta_{j,\malpha})$ to \eqref{eq_mGNPC} under   constraint \eqref{eq3} and the following addtional constraint:
\begin{equation}    \left({\malpha} \cdot \bm{q}_j = 0 \Longrightarrow {\theta}_{j,\malpha} = 0\right) \text{ and }  \left({\malpha} \cdot \bm{q}_j = K_j \Longrightarrow {\theta}_{j,\malpha} = 1\right),\label{eq10}\end{equation} 
where the additional constraint (\ref{eq10}) corresponds to the constraint  (\ref{eq10-eq}) under the GNPC setting. 
}

{Following the above discussion, both the original GNPC method and the modified GNPC method can be formulated in a unified estimation framework~\citep{Ma2023} of minimizing \eqref{eq_mGNPC} under different constraints. In particular, since $\sum_{i=1}^N \sum_{j=1}^J (R_{i,j} - \theta_{j,\malpha_i})^2=\sum_{\malpha \in \{0,1\}^K} \sum_{i: {\malpha}_i = \malpha} \sum_{j=1}^J (R_{i,j} - \theta_{j,\malpha})^2$, 
 we can rewrite   \eqref{eq_mGNPC}
  equivalently as the following loss function \begin{align}
  \ell(\A,\T|\R)=
  \sum_{i=1}^N \sum_{j=1}^J (R_{i,j} - \theta_{j,\malpha_i})^2, \label{eq8}\end{align}
where minimizing the loss function (\ref{eq8}) with respect to $({\A},{\T})$ under the constraints (\ref{eq3}) and (\ref{eq10}) obtains the original GNPC estimators  in ~\cite{chiu2019consistency} and the modified GNPC estimators in \cite{Ma2023} can be obtained by minimizing \eqref{eq8} under the constraint (\ref{eq3}) only.} 

\paragraph{Limitations of Existing  Consistency Results for Nonparametric CDMs}
Existing theoretical research has offered valuable insights into the practical utility of nonparametric methods. It has been shown that the NPC estimators are statistically consistent for estimating subjects' latent attributes under certain CDMs~\citep{wang2015}. Similarly, the GNPC estimator's ability to consistently classify subjects has been established ~\citep{chiu2019consistency}. However, current theoretical assurances for these nonparametric classification methods come with their own set of limitations. 

A fundamental assumption for the NPC method to yield a statistically consistent estimator of $\A$ is that $\P(R_{i,j}=1|\eta_{j,\malpha}=0)<0.5$ and $\P(R_{i,j}=1|\eta_{j,\malpha}=1) > 0.5$, where ${\eta_{j,\malpha}}$ represents the binary ideal responses (either 0 or 1) under the considered model~\citep{wang2015}. 
However, as previously pointed out, this binary ideal response becomes restrictive when working with more complex CDMs. The binary ideal response, limited to representing the complex latent attribute patterns of examinees through two states, could potentially oversimplify the actual complexity of the scenario. This limitation, in turn, constrains the practical application of the NPC method in instances where the underlying true model is more sophisticated. For instance,~\cite{chiu2018cognitive} provided an illustrative example highlighting this restriction, showing the possibility of misclassifications when the underlying true model is the saturated GDINA model.

Although the GNPC method addresses the oversimplification problem of the NPC method, a new restrictive assumption emerges in the existing theory for the GNPC method. {Specifically, Theorem 1 in \cite{chiu2019consistency} assumes initialization of the memberships $\hat\malpha_i^{(0)}$s that consistently estimates the ground truth in order to establish the consistency theory for GNPC. Similarly, \cite{Ma2023} assumes the existence of a calibration dataset that provides consistent estimations $\widehat{\A}_c$ for the true latent class membership $\A_c^0$ of the calibration subjects. Under these assumptions, $\widehat{\eta}_{j,\malpha}^{(w)}$ can be estimated using consistent membership estimations, which further support the consistency theory.}
The assumption concerning the existence of an initial set of consistent estimates or a calibration dataset may be restrictive and hard to satisfy in practice. To address this issue, we present new theoretical results demonstrating that the {consistency} of the GNPC method can be established without the need for a consistent initialization or a calibration dataset. These findings are detailed in the subsequent section.

\section{Main Results}\label{sec:estimator}

{Based on the unified framework of two GNPC methods outlined in Section~\ref{sec:setup}, we will establish the theoretical properties of both the original GNPC method \citep{chiu2019consistency} and the modified GNPC method \citep{Ma2023} under less stringent conditions. Regarding implementation, estimation algorithms for both the original and modified GNPC methods have been detailed in \cite{chiu2018cognitive} and \cite{Ma2023}, respectively. We recommend using these well-established methods for estimation.}
Before delving into the statistical behaviors of the aforementioned general nonparametric estimators, we outline the needed regularity conditions. Consider a model sequence indexed by \( (N,J) \), where both \( N \) and \( J \) tend to infinity, while \( K \) is held constant. For clarity, let the true parameters generating the data be represented as \( (\T^0,\Q^0,\A^0) \), and other true parameters are also denoted with superscript $0$. Assumptions are made on these true parameters as follows.

\begin{Ass}\label{Ass1}
    There exists $\delta>0$ such that \begin{align}
        \min_{1\leq j\leq J} \left\{ \min_{\malpha \circ \bm{q}_j \neq \tilde{\malpha}\circ \bm{q}_j} (\theta^0_{j,\malpha} - \theta^0_{j,\tilde{\malpha}})^2 \right\} \geq \delta .\nonumber
    \end{align}
\end{Ass}

\begin{Ass}\label{Ass2}
    There exist {$\{\delta_J:\delta_J>0\}_{J=1}^{\infty}$} and a constant $\epsilon >0 $ such that \begin{align}
        & \min_{1\leq k\leq K} \frac{1}{J} \sum_{j=1}^J \1\{\bm{q}_j^0 = \bm{e}_k\}\geq \delta_J \label{eq12};\\
        & \min_{\malpha\in\{0,1\}^K} \frac{1}{N} \sum_{i=1}^N \1\{ \malpha_i^0 = \malpha\}\geq \epsilon .\label{eq13}
    \end{align}
\end{Ass}

Assumption~\ref{Ass1} serves as an identification condition for local latent classes at each item level, ensuring that the item parameters of different local latent classes, influenced by $\K_j$, are sufficiently distinct. The gap, denoted as $\delta$, measure the separation between latent classes, thereby quantifying the strength of the signals. In the finite-$J$ regime, a $\delta>0$ is required for studying identifiability~\citep{kohn2017procedure,xu2018,gu2019b,gu2023}. 
Assumption~\ref{Ass2} pertains to the discrete structures of $\Q$ and $\A$. Here, (\ref{eq13}) implies that the $2^K$ latent patterns are not too unevenly distributed in the sample. An equivalent requirement in random-effect latent class models is $p_{\malpha}>0$ for all $\malpha \in \{0,1\}^K$, where $p_{\malpha}$ represents the population proportion of latent pattern $\malpha$. For Assumption~\ref{Ass2}, within the finite-$J$ regime, (\ref{eq12}) is similar to the requirement that ``$\Q$ should contain an identity submatrix $\mathbf{I}_K$''~\citep{chen2015, xu2018}. However, as $J$ approaches infinity, a finite number of submatrices $\mathbf{I}_K$ in $\Q$ may not be sufficient to ensure estimability and consistency. Therefore, (\ref{eq12}) necessitates that $\Q$ includes an increasing number of identity submatrices, $\mathbf{I}_K$, as $J$ grows. A similar assumption on $\Q$ was made by~\cite{wang2015} when they were establishing the consistency of the NPC method. 
{It is worth mentioning that the lower bound $\delta_J$ in \eqref{eq12} in Assumption~\ref{Ass2} is allowed to decrease to zero as $J$ goes to infinity.}

In the following subsections, we study the consistency properties of the modified GNPC method with the constraint (\ref{eq3})  and the original GNPC method with both constraints (\ref{eq3}) and (\ref{eq10}).
As the modified GNPC method involves less constraints compared to the original GNPC method, for convenience,  we first present results for the modified GNPC method  in Section \hyperlink{section3.1}{3.1} and then for the original GNPC method in Section \hyperlink{section3.2}{3.2}.

\subsection{Consistency Results for Modified GNPC}

In this section, we discuss the consistency results for the modified GNPC method. The following main theorem first validates the {consistency} of the modifed GNPC method under the constraint (\ref{eq3}) and provides a bound for its rate of convergence in recovering the latent attribute profiles. {We use $O(\cdot)$ and $o(\cdot)$ to denote the big-O and small-o notations, respectively, and $O_p(\cdot)$ and $o_p(\cdot)$ as their probability versions for convergence in probability.}

\begin{Thm}[{Consistency of modified GNPC method}] \label{Thm1}
    Consider the $(\widehat{\A},\widehat{\T}) = \argmin_{(\A,\T)}$ $\ell(\A,\T|\R)$ under the constraint (\ref{eq3}). When $N,J\to \infty$ jointly, suppose $\sqrt{J} = O\left( N^{1- c} \right)$ for some small constant $c\in(0,1)$. Under Assumption~\ref{Ass1} and Assumption~\ref{Ass2}, the {classification} error rate is \begin{align}
            \frac{1}{N} \sum_{i=1}^N \1\{ \widehat{\malpha}_i \neq \malpha_i^0\} = o_p \left( \frac{(\log J)^{\tilde{\epsilon}}}{ \delta_J\sqrt{J}}\right), \label{eq15}
        \end{align}
        \noindent where for a small positive constant $\tilde{\epsilon}>0$.
\end{Thm}

Theorem~\ref{Thm1} bounds the error of the estimator $\widehat{\A}$, which establishes the {consistency of the latent attributes} of the nonparametric method, and even allows the rate $\delta_J$ to go to zero. 
Theorem~\ref{Thm1} also offers insight into the accuracy of estimating $\A$ with finite samples and finite $J$. In particular, if $\delta_J$ is a constant, then the finite sample error bound in (\ref{eq15}) becomes $o_p((\log J)^{\tilde{\epsilon}} /\sqrt{J})$. {Ignoring the $\log$ terms, the result shows that the {classification} error rate can be dominated by the order of $J^{-1/2}$, indicating that a longer item set facilitates more accurate classification for the latent profiles of all subjects.}
Note the scaling condition that $N\exp(-J t^2)\to 0$ for any positive fixed $t>0$ in~\cite{chiu2019consistency} and~\cite{wang2015} essentially requires the growth rate of $J$ to be at least the order of $\log N$. {In contrast, our scaling condition only assumes that the number of items goes jointly with $N$ at a slower rate, which can be more easily satisfied.}

The following corollary demonstrates that under certain conditions, the item parameters can be consistently estimated via the modified GNPC method as $N,J\to\infty$.

\begin{Thm}[Item Parameters Consistency] \label{Thm2}
    Under Asssumptions~\ref{Ass1} and \ref{Ass2} and the scaling conditions given in Theorem~\ref{Thm1}, we have the the following uniform consistency result  for all $j\in[J]$ and $\malpha\in\{0,1\}^K$: \begin{align}
        \max_{j,\malpha}\abs{\hta - \theta_{j,\malpha}^0} = o_p\left(\frac{1}{\sqrt{N^{1-\tilde{c}}}} \right) + o_p\left( \frac{(\log J)^{\tilde \epsilon}}{\delta_J \sqrt{J}}\right), \nonumber
    \end{align}
    \noindent where $\tilde{c}$ and $\tilde{\epsilon}$ are  small positive constants.
\end{Thm}

This theorem builds on the {consistency result} established in Theorem~\ref{Thm1} to establish the uniform consistency in parameter estimation. The condition \(\sum_{i=1}^N \1\{\malpha_i^0 = \malpha\}\geq N\epsilon\) for all \(\malpha \in \{0,1\}^K\) ensures that there are enough samples within each class to provide accurate estimates of item parameters.  {This is reflected in the first error term \(o_p(1/\sqrt{N^{1-\tilde{c}}})\), which achieves nearly optimal  $\sqrt{N}$-consistency.} Notably, the added \(\tilde{c}\) term arises due to the number of parameters going to infinity jointly with the sample size \(N\), causing a slight deviation from the optimal error rate of \(O_p(1/\sqrt{N})\). The maximum deviation $\max_{j,\malpha}|\widehat{\theta}_{j,\malpha} - \theta_{j,\malpha}^0|$ is also influenced by the classification errors for the unknown latent attributes, which is shown in the second error term \(o_p((\log J)^{\tilde\epsilon}/(\delta_J\sqrt{J}))\). {In conclusion, the upper bound for the maximal error in estimating item parameters comprises a term that denotes nearly optimal $\sqrt{N}$-consistency, accompanied by an additional term related to the errors in classifying the latent attributes. Our theory suggests that both the sample size and the test length need to be sufficiently large to ensure accurate estimation of the item parameters, given that the latent attributes of the subjects must also be estimated.}

\subsection{Consistency Results for Original GNPC}

In this subsection, we discuss the {consistency} result of the original GNPC method. Since the original method adds an additional constraint (\ref{eq10}) compared to the modified method, which causes some of the parameters $\theta_{j,\malpha}$ to be $0$ or $1$, additional notation are needed to characterize how this potential variation affects the consistency outcome. Denote 
\begin{align}
    \lambda_{N,J}^2 = \frac{1}{NJ}\sum_{j=1}^J\left(\sum_{i: {\malpha}_i^0 \cdot \bm{q}_j = 0} (\theta_{j,{\malpha}_i^0}^0 - 0)^2  +   \sum_{i:{\malpha}_i^0 \cdot \bm{q}_j = K_j} (1-\theta_{j,{\malpha}_i^0}^0)^2\right) ,\label{eq17}
\end{align}
 which represents the average squared distance between the true parameters and the associated zero/one values. For establishing {consistency} for the original GNPC method, an additional assumption is needed.

\begin{Ass} \label{Ass3}
    For any $j\in[J]$,  we have $\theta_{j,\malpha=\mathbf 0}^0 <1/2 < \theta_{j, \malpha=\mathbf 1}^0.$
\end{Ass}
Assumption~\ref{Ass3} plays a similar role to Assumption~\ref{Ass1}, as both measure the separation between different latent classes. While this appears to be a relatively mild condition and may seem similar to the one presented in~\cite{wang2015}, as discussed in Section~\ref{sec:setup}, it remains applicable to complex CDMs. The following theorem validate the consistency of the {nonparametric} classification method under the original GNPC setting, provides a similar bound as Theorem~\ref{Thm1} for the misclassification rate.

\begin{Thm}[GNPC Consistency] \label{Thm3}
    Consider the $(\widehat{\A},\widehat{\T}) = \argmin_{(\A,\T)}\ell(\A,\T|\R)$ under the constraints (\ref{eq3}) and (\ref{eq10}). When $N,J\to \infty$ jointly, suppose $\sqrt{J} = O\left(  N^{1- c} \right)$ for some small constant $c\in(0,1)$. Under Assumptions~\ref{Ass2} and \ref{Ass3}, the {classification} error rate  is \begin{align}
            \frac{1}{N} \sum_{i=1}^N \1\{ \widehat{\malpha}_i \neq \malpha_i^0\} \leq o_p \left( \frac{(\log J)^{\tilde{\epsilon}}}{\delta_J\sqrt{J}}\right) + \frac{4\lambda_{N,J}^2}{\delta_J}. \nonumber
        \end{align}
\end{Thm}

The {classification} error rate for the original GNPC method is slightly different from the result given in Theorem~\ref{Thm1}. An extra item $4\lambda_{N,J}^2/\delta_J$ is added into the error rate. {This additional term reflects the number of items that violate constraint \eqref{eq10}, as \eqref{eq17} will be larger when there are more items with $\theta_{j,\malpha}$ that is neither 0 nor 1. The impact of the additional error introduced by $\lambda_{N,J}$ is further illustrated by Example~\ref{Eg1} and the simulation studies in Section~\ref{sec:simu}.} The details of Theorem~\ref{Thm3}'s proof can be found in Section~\ref{App2} of the supplementary material.

It is worth mentioning that without further regularity conditions, it might be challenging to avoid the additional error term $4\lambda_{N,J}^2/\delta_J$. In the existing consistency results for both the NPC and the modified GNPC methods~\citep{wang2015, chiu2019consistency}, a crucial step involves ensuring that for each examinee $i$, the true attribute profile minimizes $\E[d_i(\malpha_m)]$ across all $m$. Here, $d_i(\malpha_m) = \sum_{j=1}^J d(R_{i,j}, \widehat{\eta}_{j,\malpha_m})$ represents the distance functions used in the respective nonparametric methods. A similar approach is required in the proof of Theorem~\ref{Thm1} for the modified GNPC method. If we denote $\overline{\ell}(\A,\T) = \E[\ell(\A,\T|\R)]$ and $\overline{\ell}(\A) = \inf_{\T} \overline{\ell}(\A,\T)$, then the true latent class profiles, $\A^0$, are found to minimize $\overline{\ell}(\A)$. 

One challenge in establishing the {consistency} for the original GNPC method lies in the fact that, with the inclusion of the constraint (\ref{eq10}), the true latent class profiles $\A^0$ might not necessarily minimize $\overline{\ell}(\A)$. Let $\tilde{\A} = \argmin_{\A} \overline{\ell} (\A)$, it can be intuitively understood that $\widehat{\A}$ might approach $\tilde{\A}$ more closely than $\A^0$. Thus the additional error term originates from the discrepancy between $\overline{\ell}(\A^0)$ and $\overline{\ell}(\tilde{\A})$. Indeed, in the proof of Theorem~\ref{Thm3}, we employ the following upper bound to account for this deviation:
\begin{align}\lambda_{N,J}^2 \geq \frac{1}{NJ}\left(\overline{\ell}(\A^0)  - \overline{\ell}(\tilde{\A}) \right).\label{eq19}\end{align}

\noindent The above inequality in (\ref{eq19}) is sharp up to a constant multiple of $\lambda_{N,J}^2$, below is an illustrative example.


\begin{Eg}\label{Eg1}
    In this example, we assume that 
    the number of sample size  $N$ is $8M$ for some positive integer $M$, the number of items $J$ is 4, and the dimension of latent attribute profiles $K$ is also 4. We further assume that the four corresponding row vectors for the items in the Q-matrix are $\bmq_1 = (1,0,0,1),\ \bmq_2 = (1,1,0,0),\ \bmq_3 =(0,1,1,0)$, and $ \bmq_4 = (0,0,1,1)$, where $\bmq_j$ encodes the required latent attributes for the $j$-th item. For the true latent attribute profiles of the $8M$ samples, it is assumed that $4M$ samples exhibit latent attribute profile $(1,1,1,1)$, while the remaining $4M$ display the profile $(0,0,0,0)$. It is noteworthy that all parameters in the original GNPC method will be treated as exactly zero or one under the true latent attribute profiles $\mathbf{A}^0$ in this example, as stipulated by the constraint (\ref{eq10}). The last assumption in this example is that there exists some $ \lambda \in(0,1/2)$ such that $\theta_{j,\malpha = \bm{0}}^0 = 1/2 - \lambda$ and $  \theta_{j,\malpha = \bm{1}}^0 = 1/2 +\lambda$. Under these assumptions, the expected loss under the true latent attribute profiles satisfies \begin{align}\overline{\ell}(\A^0) - \sum_{i=1}^N\sum_{j=1}^J P_{i,j}(1-P_{i,j}) =(NJ)\left(\frac{1}{2} - \lambda\right)^2 ,\label{eq20}\end{align}   
     where $P_{i,j} := \P(R_{i,j} = 1)$ are true item response parameters, independent of the estimation process. {The derivation of \eqref{eq20} is detailed in Section~\ref{sup_sec_example1} in the Supplementary Material. } To demonstrate the sharpness of inequality (\ref{eq19}), we construct an alternative set of latent attribute profiles, denoted as $\mathbf{A}^1$. This set contains $2M$ samples of $\bm{e}_k\in\{0,1\}^4$ for each $k\in[4]$, where each $\bm{e}_k$ only contains the $k$-th latent attribute. For instance, $\bm{e}_1 = (1,0,0,0)$, $\bm{e}_2 = (0,1,0,0)$, and so on. There is a correspondence between the true latent profiles $\mathbf{A}^0$ and the constructed $\mathbf{A}^1$. Specifically, for the $2M$ samples assigned to $\bm{e}_k$ within $\mathbf{A}^1$, the true latent attribute profiles are equally divided, with half being $(0,0,0,0)$ and the other half $(1,1,1,1)$. Hence, the expected loss under the constructed latent attribute profiles fulfills \begin{align}\overline{\ell}(\A^1) - \sum_i \sum_j P_{i,j}(1- P_{i,j}) = (NJ)\cdot \left(\lambda^2+\frac{1}{8} \right).\label{eq21}\end{align}
{The derivation of \eqref{eq21} is detailed in Section~\ref{sup_sec_example1} in the Supplementary Material. }
    \noindent Thus, we have $(NJ)^{-1}\left(\overline{\ell}(\A^0) - \overline{\ell}(\A^1)\right) = -\lambda + 1/8$. Note that in this example $\lambda_{N,J}^2 = (\lambda - 1/2)^2$. If $\lambda \leq 1/13$, then one can easily verify that $-\lambda + 1/8 > \lambda_{N,J}^2/4$, and  therefore, in this case, we deduce that $$\lambda_{N,J}^2 \geq \frac{1}{NJ}\left( \overline{\ell}(\A^0) -\overline{\ell}(\tilde{\A})\right) \geq \frac{1}{NJ}\left( \overline{\ell}(\A^0) -\overline{\ell}(\A^1)\right) \geq \frac{\lambda_{N,J}^2}{4},$$
which implies the order $\lambda_{N,J}^2$ in the inequality in (\ref{eq19}) is sharp. The details of the proof can be found in Section~\ref{sup_sec_example1} of the supplementary material.
\end{Eg}

The magnitude of the additional classification error term arising from the aforementioned discrepancy is of the order $O((NJ)^{-1}(\overline{\ell}(\A^0) - \overline{\ell}(\tilde{\A})))$. As demonstrated in Example 1, $O(\lambda_{N,J}^2)$ provides a tight estimation of the order of the discrepancy $(NJ)^{-1}(\overline{\ell}( \A^0) - \overline{\ell}(\tilde{\A}))$. Therefore, the additional error term $4\lambda_{N,J}^2 /\delta_J$ in Theorem~\ref{Thm3} may not be significantly reducible.

\section{Simulation Study}\label{sec:simu}
{
In this section, we conduct a comprehensive simulation study to illustrate our theoretical findings of both the original GNPC~\citep{chiu2018cognitive} and the modified GNPC~\citep{Ma2023}. Note that in the existing literature~\citep{chiu2018cognitive,Ma2023}, various numerical studies have already demonstrated the effectiveness of both the original GNPC method and the modified GNPC method in small sample settings. Therefore, our focus here primarily lies on scenarios where both the sample size and test length are relatively large to illustrate our theoretical results. 
}

For the data-generating process, followed by the simulation design of~\cite{chiu2018cognitive} and~\cite{Ma2023}, we consider two settings: (1) items are simulated using the DINA model, and (2) items are simulated from GDINA model, as detailed in Section~\ref{sec:setup}. The manipulated conditions include: the sample size $N\in\{300, 600, 1000\}$; the test length $J \in \{50, 100, 200, 300, 400, 500\}$; the number of latent attributes $K\in\{3, 5\}$. For $K=3$, the Q-matrix is constructed with two identity $K\times K$ submatrices, and the remaining items are generated uniformly from all possible non-zero patterns. It is worth mentioning that this generating process adheres to \eqref{eq12} in Assumption~\ref{Ass2}. For the case of $K=5$, the Q-matrix is restricted to contained items that measure up to 3 attributes and constructed the same way as that for $K=3$. For the data conforming to the DINA model, we simulate $s_j$ and $g_j$ independently from a uniform distribution $Unif[0,r]$ with $r\in\{0.2,0.4\}$. {For data conforming to the GDINA model, the item parameters are simulated following the framework outlined in \cite{chiu2018cognitive} as follows. For any item $j$, let $K_j^* = \sum_{k=1}^Kq_{jk}$ be the number of required attributes of item $j$, where $q_{jk}$ is the $(j,k)$th entry in the Q-matrix. Without loss of generality, we assume that these attributes with $q_{jk}=1$ are the first $K^*_j$ attributes. For instance, if $K^*_j = 3$, i.e., item $j$ requires three attributes, we then denote the possible proficiency classes as $\malpha_1^*=(000)$, $\malpha_2^*=(100)$, $\malpha_3^*=(010)$, $\malpha_4^*=(110)$, $\malpha_5^*=(001)$, $\malpha_6^*=(101)$, $\malpha_7^*=(011)$, and $\malpha_8^*=(111)$. The item parameters for item $j$ are specified by the probabilities of making the correct responses for all $\malpha_i^*$ with $1\le i\le 8$. If $K_j^*=2$, we only need to specify the probabilities for $\malpha_i^*$ with $1\le i\le 4$ since the remaining attributes are irrelevant for distinguishing among the proficiency classes, and if $K_j^*=1$, we only need to specify the probabilities for $\malpha_1^*$ and $\malpha_2^*$ \citep{chiu2018cognitive}. Analogous to the data generation process under the DINA model, we simulate two noise levels under the GDINA model as in \cite{Ma2023}, with item parameters provided in Table~\ref{table_small_noise} for small noises and Table~\ref{table_large_noise} for large noises, respectively. Note that Table~\ref{table_small_noise} contains seven rows, while Table~\ref{table_large_noise} contains six, each row representing a distinct set of item parameters. For each noise level, the set of item parameters for each item $j$ is sampled randomly from those rows with $K^*=K_j^*$ in each table.}

\begin{table}[ht]
\centering
\begin{tabular}{ccccccccc}
\hline
& $P(\malpha_1^*)$ & $P(\malpha_2^*)$ & $P(\malpha_3^*)$ & $P(\malpha_4^*)$ & $P(\malpha_5^*)$ & $P(\malpha_6^*)$ & $P(\malpha_7^*)$ & $P(\malpha_8^*)$ \\ \hline
&0.2 & 0.9 &  &  &  & & & \\
  $K^*=1$ &0.1 & 0.8 &  & & & & & \\
&0.1 & 0.9 &  &  &  & & & \\\hdashline
&0.2 & 0.5 & 0.4 & 0.9 & & & & \\
   $K^*=2$  &0.1 & 0.3 & 0.5 & 0.9 & & & & \\
&0.1 & 0.2 & 0.6 & 0.8 & & & & \\\hdashline
   $K^*=3$  &0.1 & 0.2 & 0.3 & 0.4 & 0.4 & 0.5 & 0.7 & 0.9 \\ \hline
\end{tabular}
\caption{Item response parameters for GDINA with small noises, where $K^*$ donotes the number of required attributes of a considered item.}\label{table_small_noise}
\end{table}

\begin{table}[ht]
\centering
\begin{tabular}{ccccccccc}
\hline
&$P(\malpha_1^*)$ & $P(\malpha_2^*)$ & $P(\malpha_3^*)$ & $P(\malpha_4^*)$ & $P(\malpha_5^*)$ & $P(\malpha_6^*)$ & $P(\malpha_7^*)$ & $P(\malpha_8^*)$ \\ \hline
\multirow{2}{*}{   $K^*=1$ }&0.3 & 0.7 &  & & & & & \\
&0.3 & 0.8 &  &  &  & & & \\\hdashline
&0.3 & 0.4 & 0.7 & 0.8 & & & & \\
   $K^*=2$ &0.3 & 0.4 & 0.6 & 0.7 & & & & \\
&0.2 & 0.3 & 0.6 & 0.7 & & & & \\\hdashline
   $K^*=3$ &0.2 & 0.3 & 0.3 & 0.4 & 0.4 & 0.5 & 0.6 & 0.7 \\ \hline
\end{tabular}
\caption{Item response parameters for GDINA with large noises, where $K^*$ donotes the number of required attributes of a considered item.}\label{table_large_noise}
\end{table}

{For the latent attribute patterns, they are generated using either a uniform setting, where each proficiency class is drawn with a uniform probability of $2^{-K}$, or a multivariate normal threshold model as described by~\cite{chiu2018cognitive}. In this model, each subject's attribute profile is linked to a latent continuous ability vector $\bm{z} \sim \mathcal{N}_K(\bm{0}, \bm{\Sigma})$. The diagonal elements of $\bm{\Sigma}$ are fixed at 1, while the off-diagonal elements are set to $0.3$ for a low-correlation scenario and $0.7$ for a high-correlation scenario. The attribute profile is then derived from $\bm{z}$ by applying a truncation process as follows:
\begin{equation*}
    \alpha_{ik}=\left\{\begin{aligned}
        &1,\quad z_{ik}\ge \Phi^{-1}\Big(\frac{k}{K+1}\Big);
        \\&0,\quad \text{otherwise.}
    \end{aligned}\right.
\end{equation*}
where $\Phi$ is the cumulative distribution function of the standard normal distribution.}

To illustrate our theoretical results, we compute the pattern-wise agreement rate (PAR):
\begin{equation*}
    \text{PAR}=\frac1{N}{\sum_{i=1}^NI\{\hat\malpha_i=\malpha_i\}}.
\end{equation*} We apply the original GNPC  and modified GNPC estimation methods under all manipulated scenarios. {Following the algorithms proposed by \cite{chiu2018cognitive} and \cite{Ma2023}, both methods are initialized using latent attributes estimated with the NPC method for computational efficiency}. In each scenario, we conduct 100 replications and calculate the mean value of the PARs. The iteration process is terminated when
$ 
    \frac1{N}{\sum_{i=1}^NI\{\hat\malpha_i^{(t-1)}\neq \hat\malpha_i^{(t)}\}}<0.001
$  or exceeding the maximal number of iterations set as 500. {Additionally, we conducted simulations in these scenarios where both methods are initialized with random latent attributes, resulting in estimation errors similar to those obtained with NPC initialization. Details and results are provided in Section~\ref{supp_sec_simu} of the Supplementary Material.}

\begin{figure}[htbp!]
\centering    
\subfigure{
        \includegraphics[width=6in,height=1.6in]{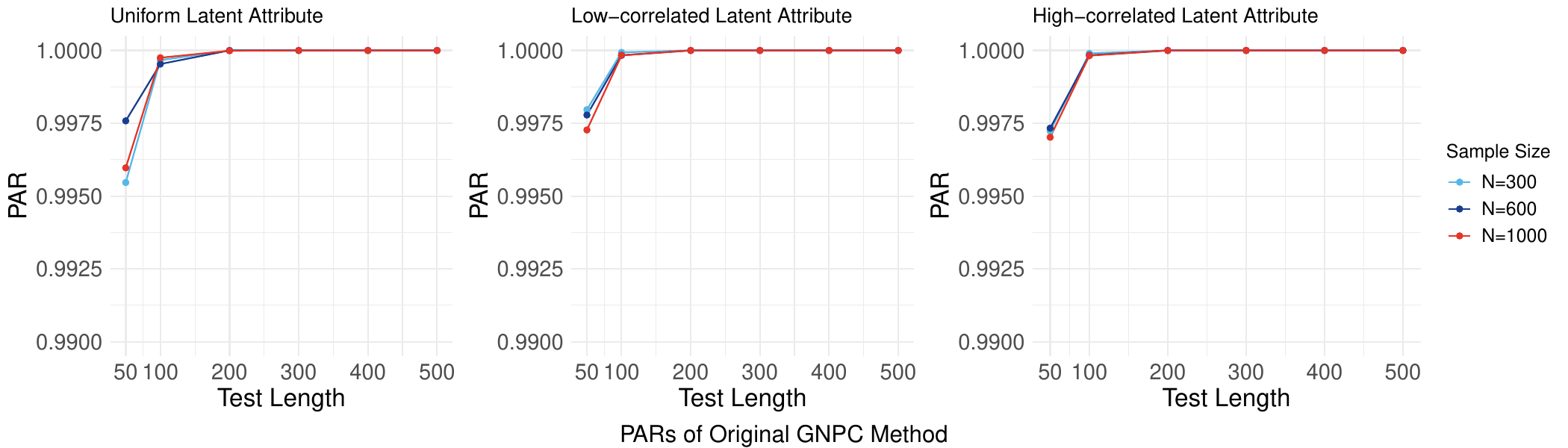}}\hspace{0.2in}
        
        \subfigure{
        \includegraphics[width=6in,height=1.6in]{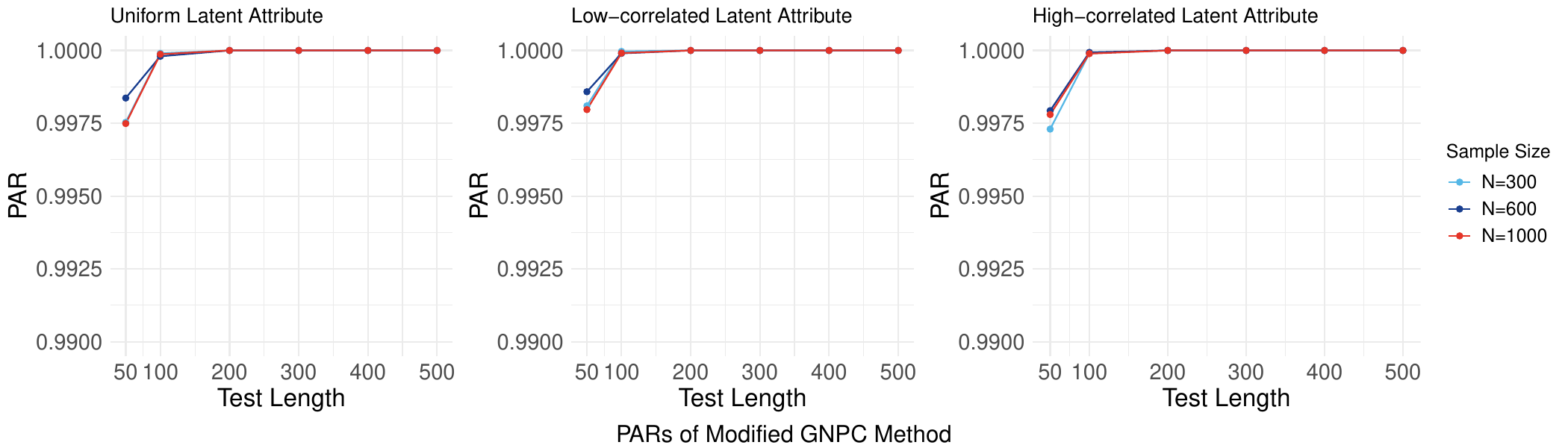}}\hspace{0.2in}
    \caption{{PARs when the data are generated using the DINA model with $K=3$ and $r=0.2$.}}
    \label{fig:K3r0.1}
\end{figure}

\begin{figure}[htbp!]
\centering    
\subfigure{
        \includegraphics[width=6in,height=1.6in]{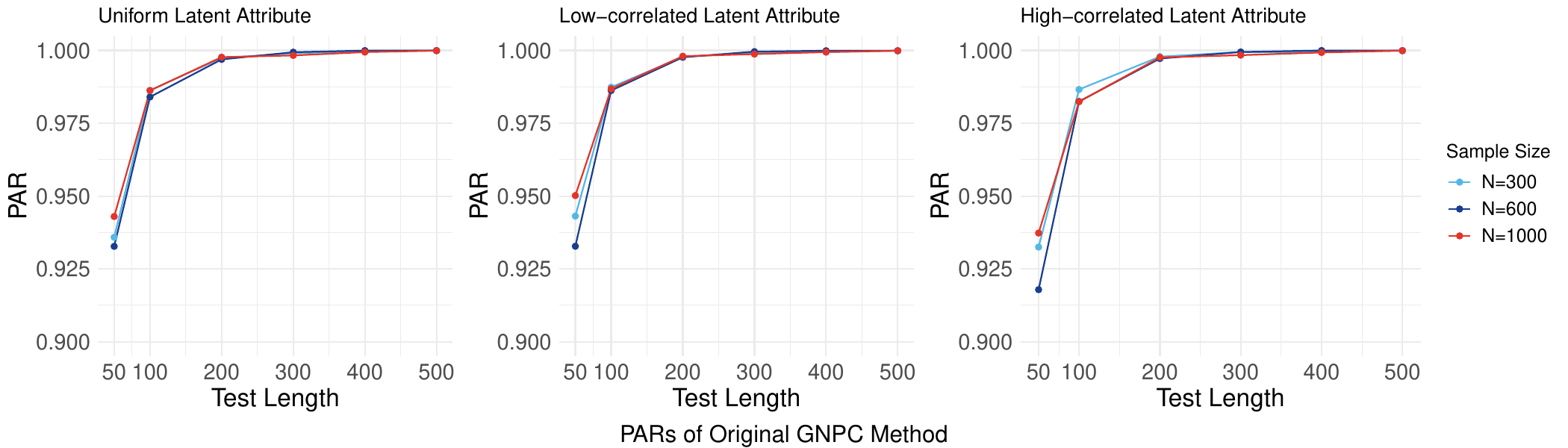}}\hspace{0.2in}
        
        \subfigure{
        \includegraphics[width=6in,height=1.6in]{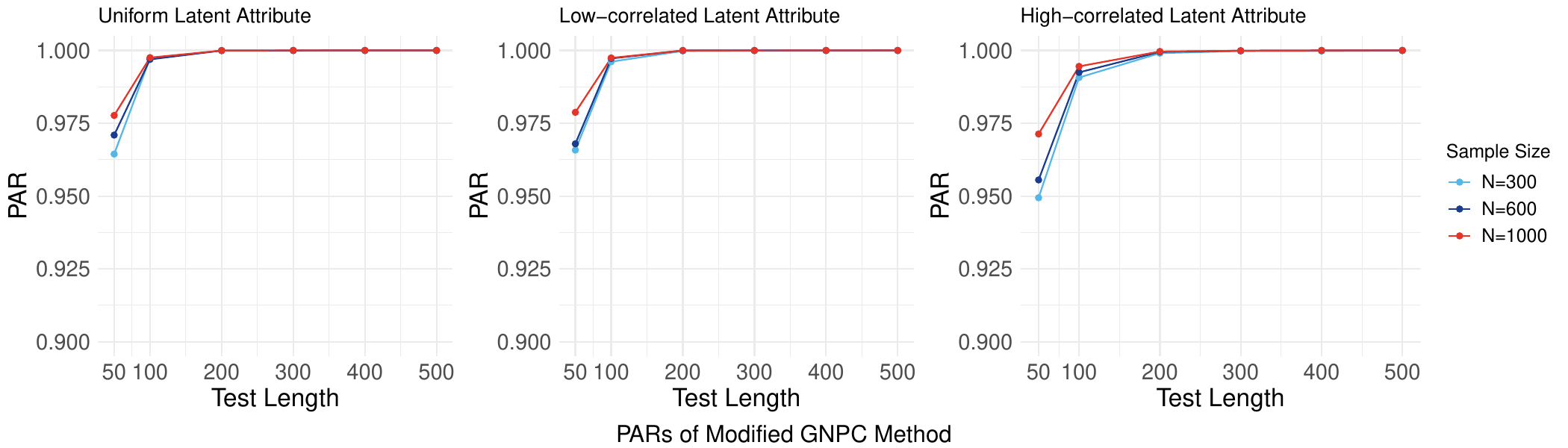}}\hspace{0.2in}
    \caption{{PARs when the data are generated using the DINA model with $K=3$ and $r=0.4$.}}
    \label{fig:K3r0.3}
\end{figure}

\begin{figure}[htbp!]
\centering    
\subfigure{
        \includegraphics[width=6in,height=1.6in]{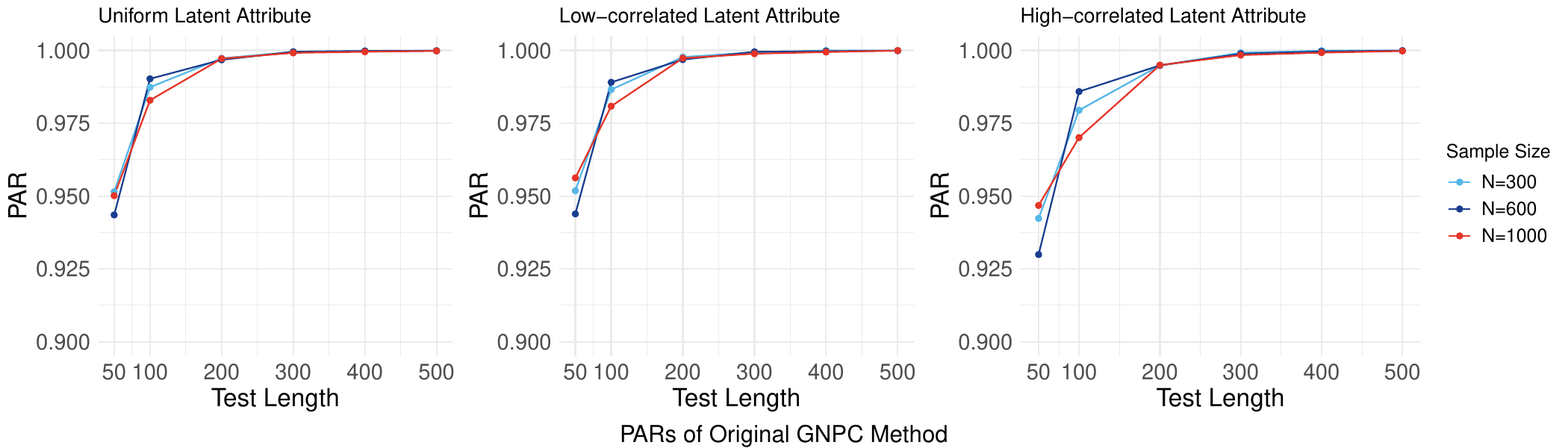}}\hspace{0.2in}
        
        \subfigure{
        \includegraphics[width=6in,height=1.6in]{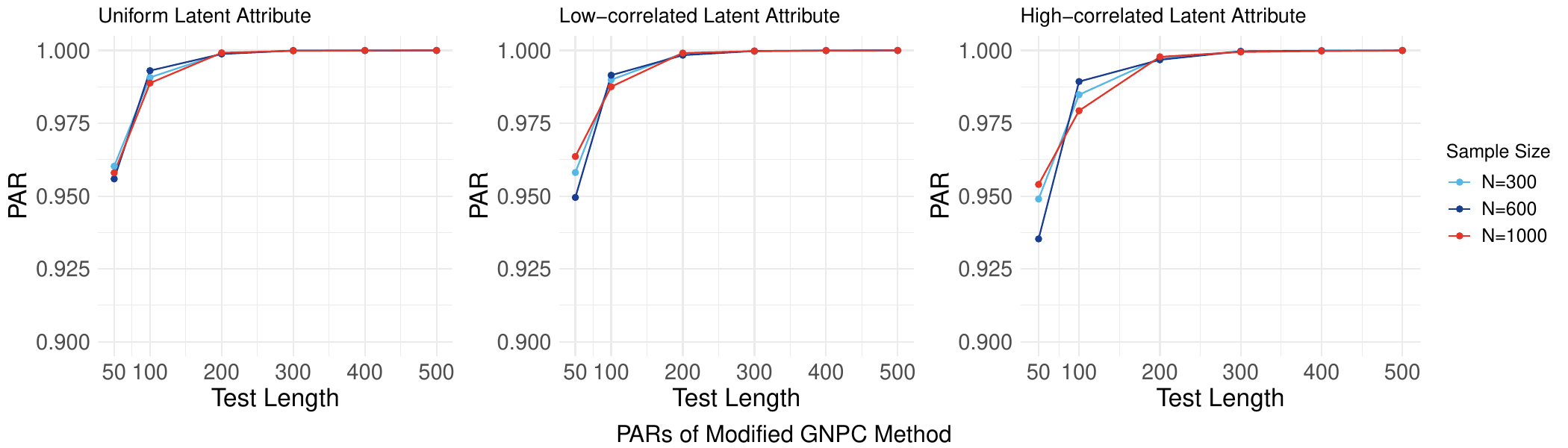}}\hspace{0.2in}
    \caption{{PARs when the data are generated using the DINA model with $K=5$ and $r=0.2$.}}
    \label{fig:K5r0.1}
\end{figure}

\begin{figure}[htbp!]
\centering    
\subfigure{
        \includegraphics[width=6in,height=1.6in]{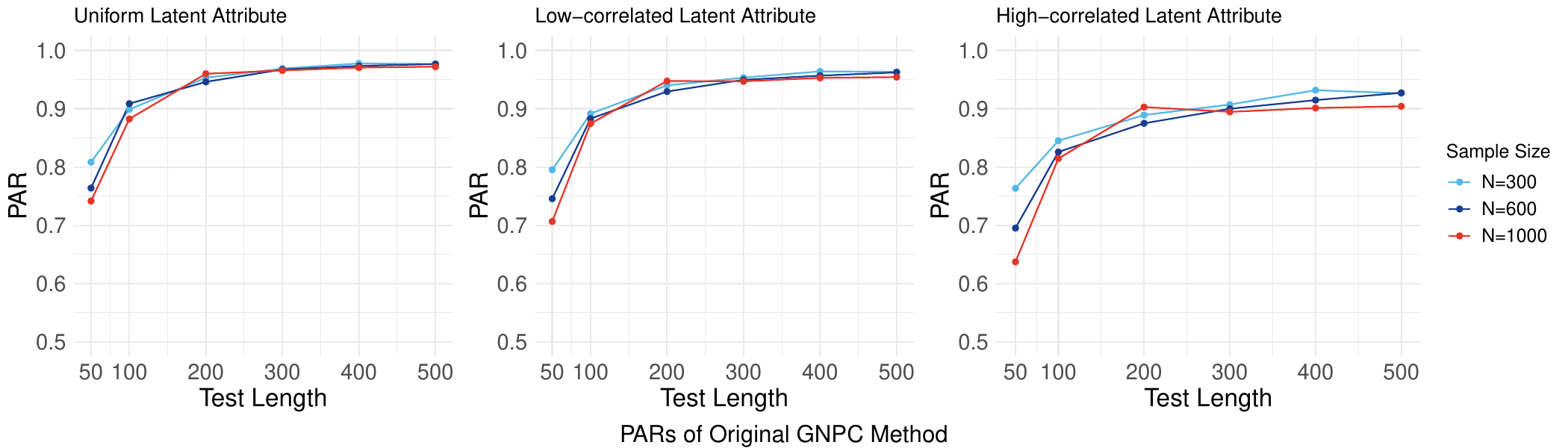}}\hspace{0.2in}
        
        \subfigure{
        \includegraphics[width=6in,height=1.6in]{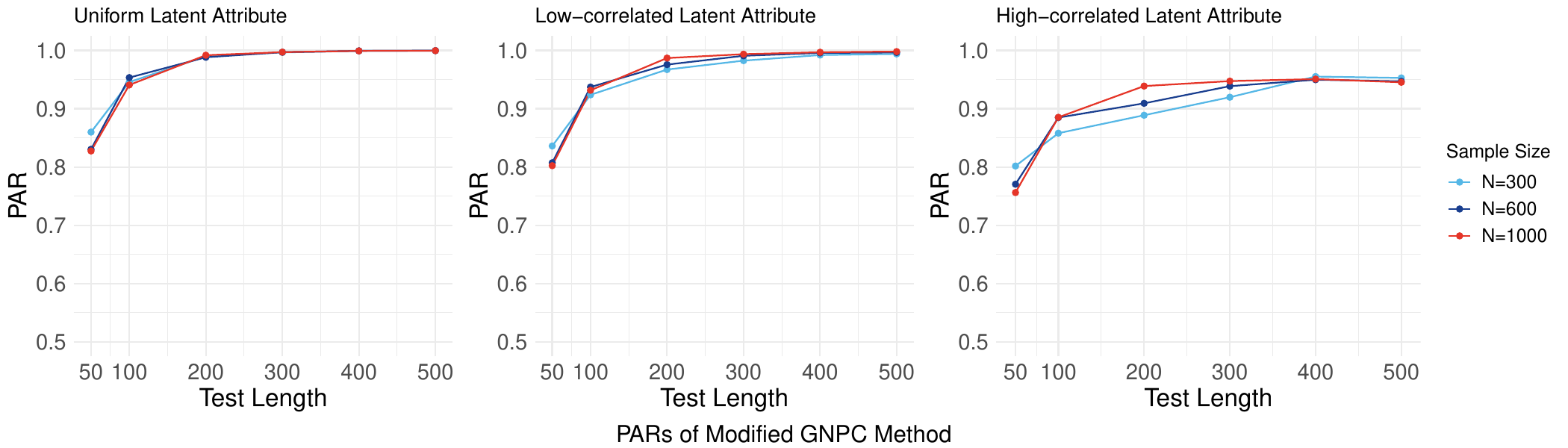}}\hspace{0.2in}
    \caption{{PARs when the data are generated using the DINA model with $K=5$ and $r=0.4$.}}
    \label{fig:K5r0.3}
\end{figure}

\begin{figure}[htbp!]
\centering    
\subfigure{
        \includegraphics[width=6in,height=1.6in]{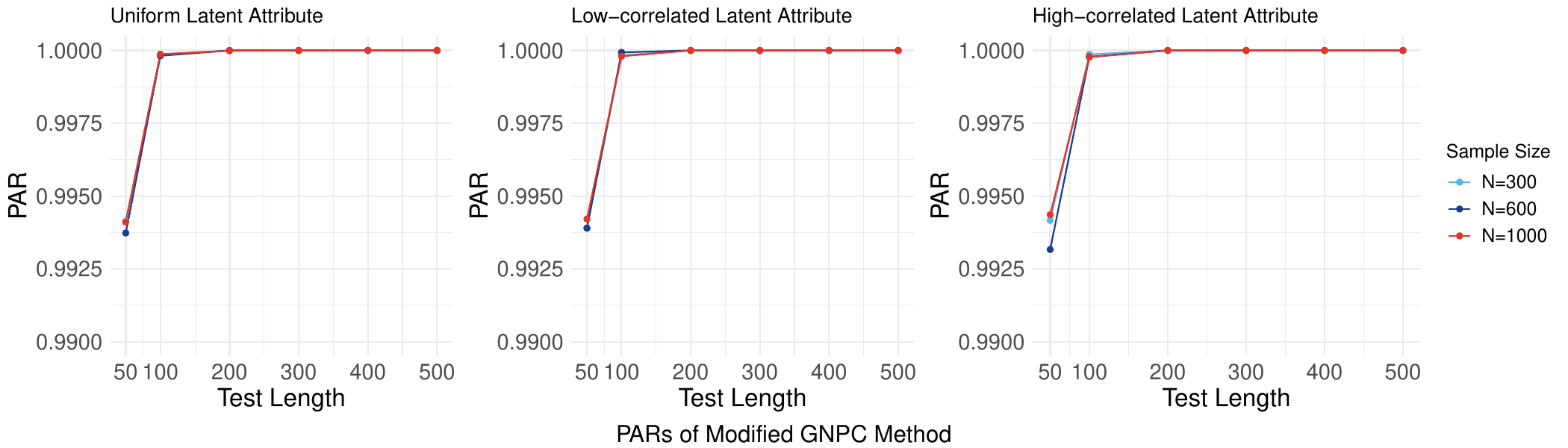}}\hspace{0.2in}

        \subfigure{
        \includegraphics[width=6in,height=1.6in]{Gk3r1mGNPC.pdf}}\hspace{0.2in}
    \caption{{PARs when the data are generated using the GDINA model with small noises and $K=3$.}}
    \label{fig:gK3r0.1}
\end{figure}

\begin{figure}[htbp!]
\centering    
\subfigure{
        \includegraphics[width=6in,height=1.6in]{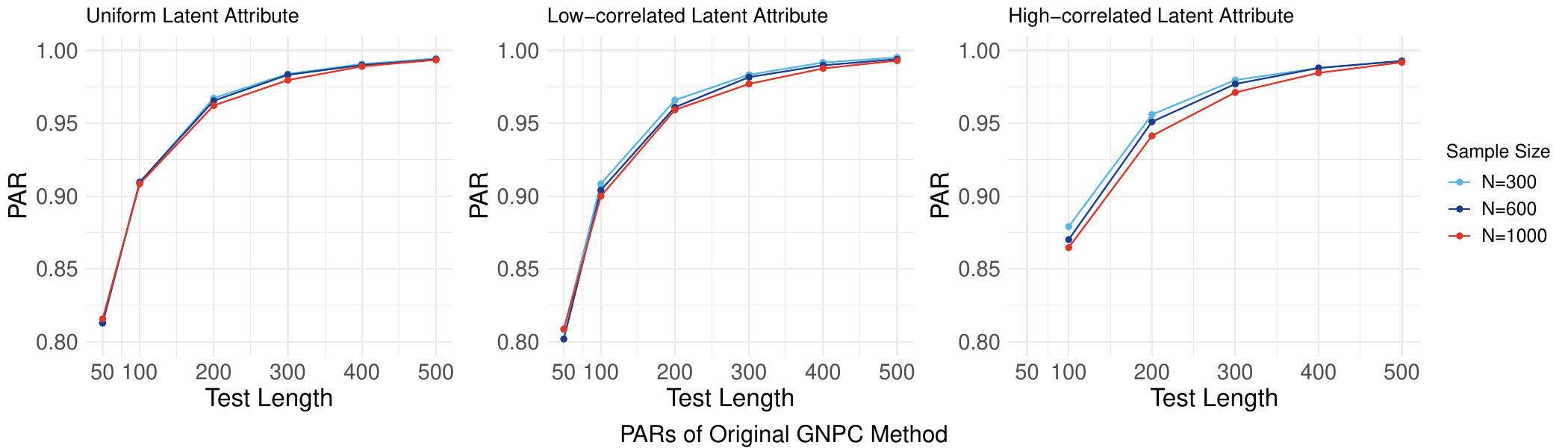}}\hspace{0.2in}
        
        \subfigure{
        \includegraphics[width=6in,height=1.6in]{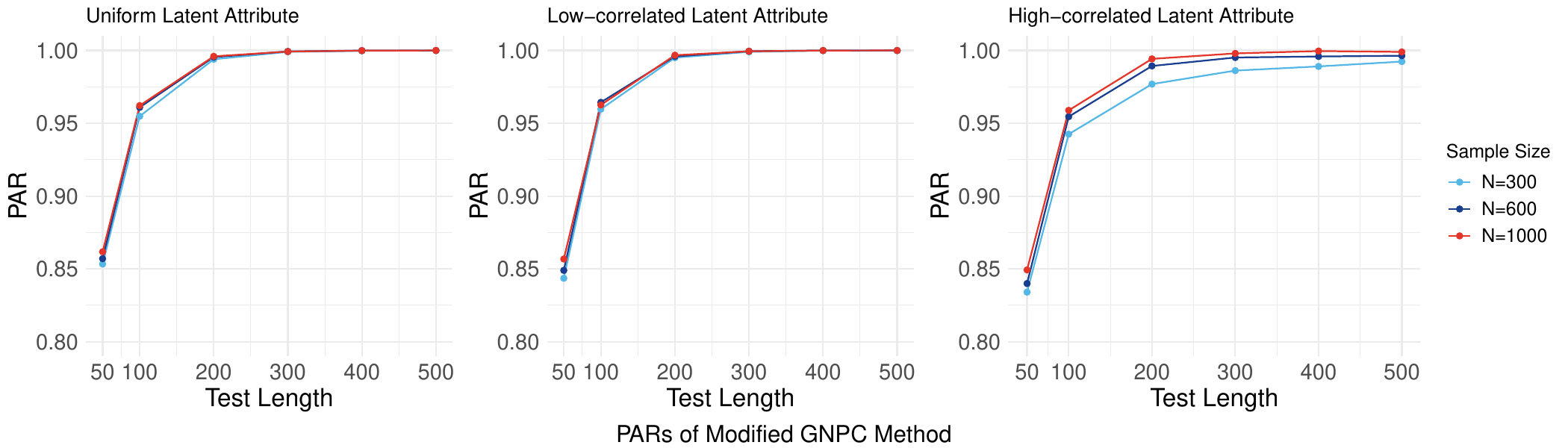}}\hspace{0.2in}
    \caption{{PARs when the data are generated using the GDINA model with large noises and $K=3$.}}
    \label{fig:gK3r0.3}
\end{figure}

\begin{figure}[htbp!]
\centering    
\subfigure{
        \includegraphics[width=6in,height=1.6in]{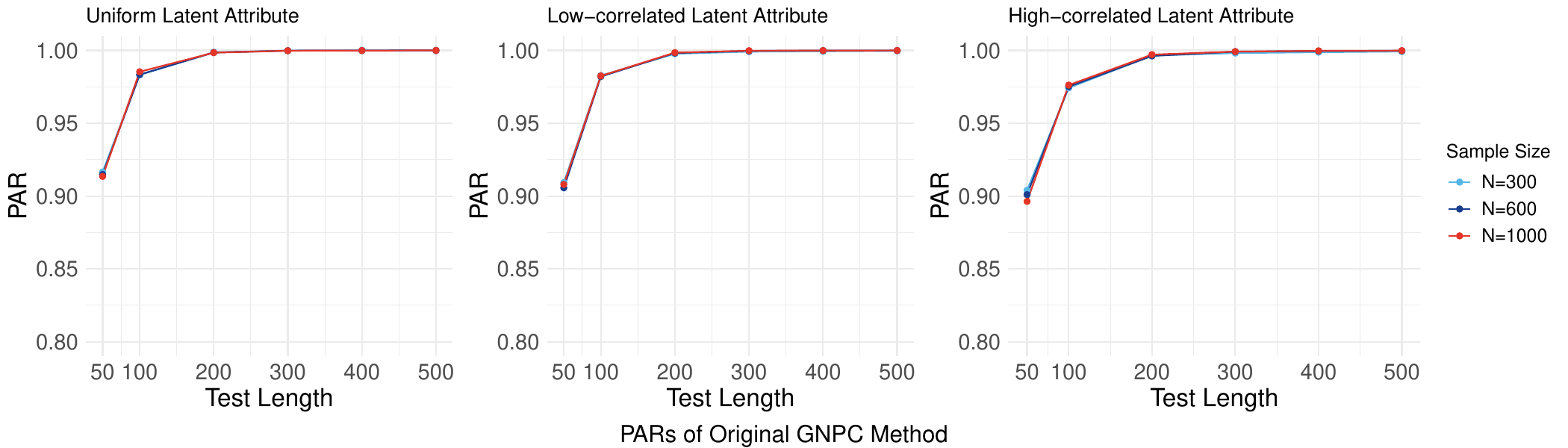}}\hspace{0.2in}
        
        \subfigure{
        \includegraphics[width=6in,height=1.6in]{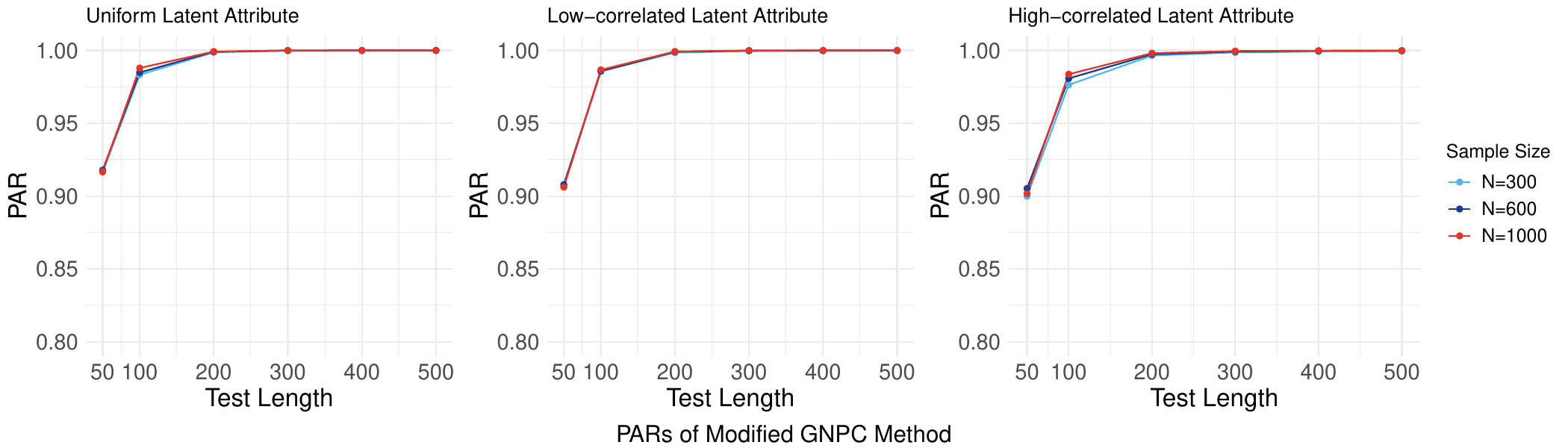}}\hspace{0.2in}
    \caption{{PARs when the data are generated using the GDINA model with small noises and $K=5$.}}
    \label{fig:gK5r0.1}
\end{figure}

\begin{figure}[htbp!]
\centering    
\subfigure{
        \includegraphics[width=6in,height=1.6in]{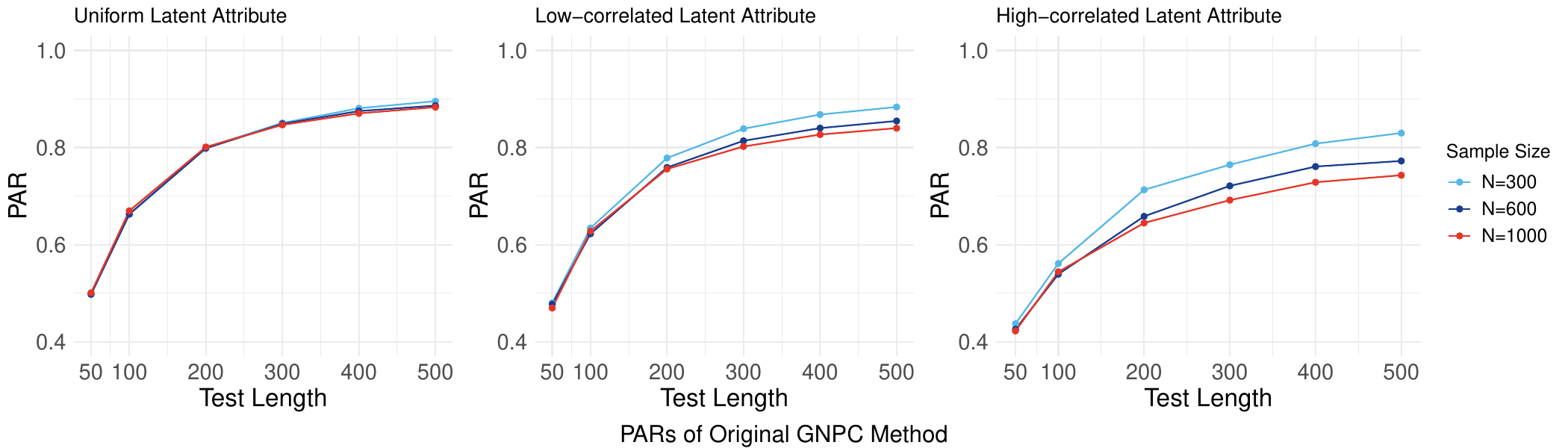}}\hspace{0.2in}
        
        \subfigure{
        \includegraphics[width=6in,height=1.6in]{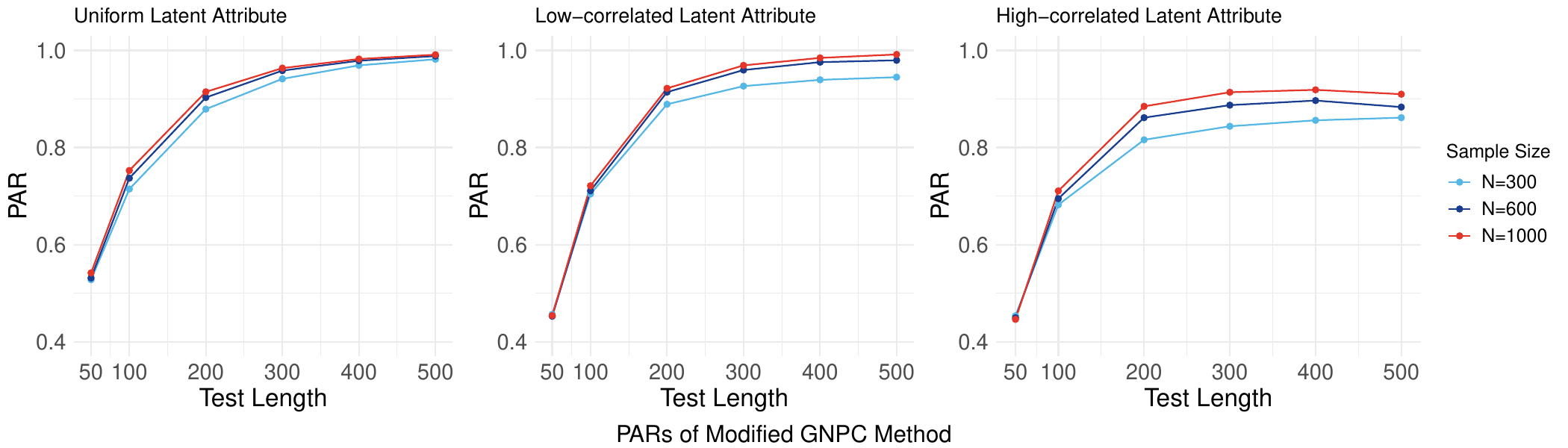}}\hspace{0.2in}
    \caption{{PARs when the data are generated using the GDINA model with large noises and $K=5$.}}
    \label{fig:gK5r0.3}
\end{figure}

{Figures~\ref{fig:K3r0.1}--\ref{fig:K5r0.3} present the PAR results when the data are generated under the DINA model, and the PAR results under the GDINA model are shown in Figures~\ref{fig:gK3r0.1}--\ref{fig:gK5r0.3}. In each figure, the upper panel presents the estimation results using the original GNPC method, while the lower panel displays the results using the modified GNPC method. From left to right, the subfigures in each row illustrate the estimations for latent attributes, simulated under three different settings: uniform, low correlation setting, and high correlation setting.}

 In general, both the original and modified GNPC methods perform well under all the model settings. The modified GNPC method exhibits a slight edge in more complex scenarios, as demonstrated in Figure~\ref{fig:K5r0.3} and Figure~\ref{fig:gK5r0.3}. A consistent trend across all figures is that as the test length $J$ increases, the PARs improve, supporting our theory that the upper bound for classification error decreases with $J$. When the data are simulated from the GDINA models, there is a slight increase in classification errors for both methods compared to those generated using the DINA models. 
In addition, comparisons between figures with lower noise levels (Figures~\ref{fig:K3r0.1}, \ref{fig:K5r0.1}, \ref{fig:gK3r0.1} and \ref{fig:gK5r0.1}) and those with higher ones (Figures~\ref{fig:K3r0.3}, \ref{fig:K5r0.3}, \ref{fig:gK3r0.3} and \ref{fig:gK5r0.3}) reveal lower classification errors with decreased noise.  In particular,  Figures~\ref{fig:K3r0.1} and \ref{fig:gK3r0.1} show nearly perfect classification results under low noise and $K=3$ settings.
Moreover, increasing the number of latent attributes typically results in less precise estimation, as evidenced by the comparisons between the settings of $K=3$ (Figures~\ref{fig:K3r0.1}, \ref{fig:K3r0.3}, \ref{fig:gK3r0.1} and \ref{fig:gK3r0.3}) and $K=5$ (Figures~\ref{fig:K5r0.1}, \ref{fig:K5r0.3}, \ref{fig:gK5r0.1} and \ref{fig:gK5r0.3}). 
Within each figure, a slight decrease in PARs is observed when the latent attributes exhibit a higher correlation. 
{When the data are simulated under larger noises and more attributes (Figure~\ref{fig:K5r0.3} and \ref{fig:gK5r0.3}), PARs from the original GNPC method appear not converging to 1 even when the sample size is $1000$ and test length is $500$. This is likely attributable to the additional error term related to $\lambda_{N,J}/\delta_J$ in Theorem~\ref{Thm3}. Notably, $\lambda_{N,J}$ in \eqref{eq17} can become large when a significant proportion of $\theta_{j,\malpha_i^0}$ fails to satisfy the constraint (\ref{eq10}).}



\section{Discussion}\label{sec:discussion}

In this work, we revise the {consistency} results for the General Nonparametric Classification (GNPC) method, originally offered in~\cite{chiu2019consistency}, under relaxed and more practical assumptions. We deliver finite sample error bounds for the considered two versions of the GNPC method. These bounds not only guarantee asymptotic consistency in estimating the latent profiles of subjects but also offer insights into the precision of these estimates in small sample situations. Furthermore, we derive uniform convergence of item response parameters $\widehat{\T}$ for the modified GNPC method. Notably, all of these advancements are achieved without the requirement for a calibration dataset.

The findings in this study open up several possibilities for future exploration. Using the consistency and finite sample error bounds established for estimating the discrete latent structure $\A$, future work can examine statistical inference on CDMs with a large number of test items and latent attributes. 
Additionally, it is important to note that in practical situations, the Q-matrix may not always be readily available. Various estimation techniques have been proposed in the literature~\citep{liu2012data, chen2015, chen2018bayesian, xu2018,li2022learning,gu2023,ma2023learning,kohn2024two}. 
This leads to a potential future direction of developing theories and computational methods for CDMs estimation with an unknown Q-matrix within the nonparametric framework.

\bibliographystyle{apalike}
\bibliography{reference}

\newpage

\appendix
{\centering \large \bf Supplementary Material for ``Clustering Consistency of General Nonparametric Classification Methods in Cognitive Diagnosis''}

\doublespacing

\bigskip

\renewcommand{\theequation}{\Alph{section}.\arabic{equation}}
\hypertarget{Appendix}{}

In this supplementary material, we provide   proofs for the theoretical results of our study. Section~\ref{sup_sec_prelim} provides notations that are used throughout the Supplementary Material. The proofs for Theorem~\ref{Thm1} and Theorem~\ref{Thm2} are presented in Section~\ref{App1}.
Section~\ref{App2} presents the proofs for Theorem~\ref{Thm3} and Section~\ref{sup_sec_example1} provides details of derivations in Example \ref{Eg1}. Section~\ref{supp_sec_simu} presents additional simulation results.

\section{Preliminaries}\label{sup_sec_prelim}

Motivated by the constraint (\ref{eq3}), we introduce the concept of a ``local'' latent class at the item level. Considering item $j$ with $Q$-vector $\bmq_j$, the constraint (\ref{eq3}) divides the collection of latent attribute profiles $\malpha$, which is $\{0,1\}^K$, based on an equivalence relationship where $\malpha \sim_j \tilde{\malpha}$ is defined by $\malpha \circ \bmq_j = \tilde{\malpha} \circ \bmq_j$, here the subscript $\sim_j$ emphasizes that the equivalence relationship is determined by the $j$-th item $\bmq_j$. On this basis, we introduce a function $\xi: \{0,1\}^K \times \{0,1\}^K \to \mathbb{N}$ where $\xi(\bm{q}_j , \malpha) = \xi(\bm{q}_j,\tilde{\malpha}) $ is equivalent to $ \malpha \circ \bm{q}_j = \tilde{\malpha} \circ \bm{q}_j$. This function assigns numbers to these equivalence classes induced by item $j$ based on some specific rules. In the following context, we call $\xi(\bmq_j^0,\malpha)$ the local latent class of $\malpha$ induced by item $j$. 
It is straightforward to verify that the number of the local latent classes induced by item $j$, denoted by $|\xi(\bmq_j,\{0,1\}^K)|$, equals $L_j = 2^{K_j}$. Here, $K_j = \sum_{k=1}^K q^0_{j,k}$ represents the number of the required latent attributes for item $j$. Consequently, we let the range of the function $\xi$ satisfies $\xi(\bmq_j,\{0,1\}^K) = [L_j] := \{1,\dots,L_j\}$. Since the local latent classes are identified up to permutations on $[L_j]$ due to their categorical nature, the mapping rules between $\xi(\bmq_j, \{0,1\}^K)$ and $[L_j]$ do not need to be completely specified in the discussion of the modified GNPC method. However, as indicated by the constraint (\ref{eq10}) and Assumption~\ref{Ass3}, the two latent attribute profiles $(0,\dots,0)$ and $(1,\dots,1)$ induce the deviation between the modified GNPC and the original GNPC methods. Thus, we specify $\xi(\bmq_j, \bm{0})=1$ and $\xi(\bmq_j, \bm{1})=L_j$ for any item $j$ in the discussion of the original GNPC method.

For brevity, we use a general notation $\mathbf{Z}=(z_{i,j})$ to denote the collection of the local latent classes for all items $j \in [J]$ and subjects $i \in [N]$, where $z_{i,j}$ represents $\xi(\bmq_j^0,\malpha_i)$. Given that $\xi(\bmq_j^0,\malpha) = \xi(\bmq_j^0,\tilde{\malpha})$ implies $\theta_{j,\malpha} = \theta_{j,\tilde{\malpha}}$ by the definition of $\xi$, we express $\theta_{j,\malpha_i}$ as $\theta_{j,z_{i,j}}$ to directly incorporate the constraint (\ref{eq3}) into the loss function (\ref{eq8}). For further notational simplicity, we may sometimes write $\theta_{j,z_{i,j}}$ simply as $\theta_{j,z_i}$. Consequently, we define
\begin{align}
    P_{i,j} = \P(R_{i,j}=1) = \theta^0_{j,z_i^0}. \tag{A.1}
\end{align}
Then the loss function (\ref{eq8}) can be rewritten as \hypertarget{A.2}{}\begin{align}
    \ell(\A,\T|\R) = \sum_{i=1}^N \sum_{j=1}^J (R_{i,j} - \theta_{j,z_i})^2. \tag{A.2}
\end{align}

\noindent Observe that $R_{i,j}^2 = R_{i,j}$, and $\E[R_{i,j}]=P_{i,j}$, we denote the expectation of the above $\ell(\A,\T|\R)$ by \hypertarget{A.3}{}\begin{align}
    \overline{\ell}(\A,\T) := \E[\ell(\A,\T|\R)] = \sum_{i=1}^N \sum_{j=1}^J (P_{i,j} - \theta_{j,z_i})^2 + \sum_{i=1}^N \sum_{j=1}^J P_{i,j}(1-P_{i,j}). \tag{A.3}
\end{align}

Note $\mathbf{Z}=(z_{i,j})$ is only determined by $\mathbf{A}$ since $\Q^0$ is known. In the subsequent context, the quantities that are determined by the latent attribute profiles $\mathbf{A}$ are sometimes denoted with a superscript $\mathbf{A}$ to emphasize their relationships with $\mathbf{A}$. Given an arbitrary $\A$, denote \begin{align}
    \ell(\A) &= \inf_{\T} \ell(\A,\T|\R) = \ell(\A,\widehat{\T}^{(\A)}|\R);  \tag{A.4} \\
    \overline{\ell}(\A) &= \inf_{\T} \overline{\ell}(\A,\T) = \overline{\ell}(\A,\overline{\T}^{(\A)}), \tag{A.5}
\end{align}

\noindent where $\widehat{\T}^{(\A)} := \argmin_{\T} \ell(\A,\T|\R)$ and $\overline{\T}^{(\A)} := \argmin_{\T}\overline{\ell}(\A,\T)$. Then under any realization of $\mathbf{A}$, the following equations hold for any local latent class $a\in[L_j]$: \hypertarget{A.6}{}\begin{align}
\widehat{\theta}_{j,a}^{(\A)} = \frac{\sum_{i=1}^N \1\{z_{i,j}^{(\A)} = a\}R_{i,j}}{\sum_{i=1}^N \1\{z_{i,j}^{(\A)}=a\}},\quad \overline{\theta}_{j,a}^{(\A)} = \frac{\sum_{i=1}^N \1\{z_{i,j}^{(\A)} = a\}P_{i,j}}{\sum_{i=1}^N \1\{z_{i,j}^{(\A)}=a\}}. \tag{A.6}
\end{align}

\noindent To derive (\hyperlink{A.6}{A.6}), note the sum $\sum_{j=1}^J\sum_{i=1}^N (R_{i,j} - \theta_{j,z_i})^2 $ equals the sum $ \sum_{j=1}^J \sum_{a=1}^{L_j} \sum_{z_i=a} (R_{i,j} - \theta_{j,a})^2$. When estimating $\widehat{\theta}_{j,a}$, we focus on minimizing the term $\sum_{z_i=a} (R_{i,j} - \theta_{j,a})^2$. For $\overline{\theta}_{j,a}$, note that $\sum_i\sum_j P_{i,j}(1-P_{i,j})$ is independent of the estimation processes. A useful observation is that, by plugging $\E[R_{i,j}]=P_{i,j}$ into (\hyperlink{A.6}{A.6}), we can find that $\E[\hta] = \bta$ holds for any $(j,a)$.

\section{Proofs of Theorem~\ref{Thm1} and Theorem~\ref{Thm2}}\label{App1}

We first outline the main steps of the proof of Theorem~\ref{Thm1} as follows and then proceed one by one.

\subsection{Outline of the first half of the proof}

\textbf{Step 1}: Express $\overline{\ell}(\A) - \ell(\A)$ by $\sum_{j=1}^J \sum_{a=1}^{L_j} n_{j,a} (\hta - \bta)^2 + \E[X]-X$, where $X := \sum_i\sum_j R_{i,j}(1-2\overline{\theta}_{j,z_i})$ depending on $\R$ and $\overline{\T}^{(\A)}$ under $\A$, and $n_{j,a} := \sum_{i=1}^N \1\{ z_{i,j}^{(\A)} = a  \}$.

\vspace{5pt}
\noindent \textbf{Step 2}: Bound $\sum_j\sum_a n_{j,a}(\hta-\bta)^2$ and $\abs{X-\E[X]}$ separately to obtain a uniform convergence rate $\sup_{\A}\abs{\overline{\ell}(\A) - \ell(\A)} = o_p(\delta_{N,J})$

\vspace{5pt}
\noindent \textbf{Step 3}: By noting the closed form of $\overline{\ell}(\mathbf{A}) - \overline{\ell}(\mathbf{A}^0) = \sum_i \sum_j (P_{i,j} - \overline{\theta}_{j,z_i}^{(\mathbf{A})})^2$ for any $\mathbf{A}$, we deduce that $\overline{\ell}(\mathbf{A}) \geq \overline{\ell}(\mathbf{A}^0)$ holds for all $\mathbf{A}$. Based on the definition of $\widehat{\mathbf{A}}$, it follows that $0 \leq \overline{\ell}(\widehat{\mathbf{A}}) - \overline{\ell}(\mathbf{A}^0) \leq 2\sup_{\mathbf{A}} |\overline{\ell}(\mathbf{A}) - \ell(\mathbf{A})| = o_p(\delta_{N,J})$, which controls the deviation $\overline{\ell}(\widehat{\mathbf{A}}) - \overline{\ell}(\mathbf{A}^0)$.

\vspace{5pt}

In some classical statistical inference contexts, consistency results for the parameters of interest are typically established through the uniform convergence of random functions associated with these parameters. For instance, if $\sup_{\vtheta \in \T} \abs{\widehat{\ell}(\vtheta) - \ell(\vtheta)}\pto 0$, and if we further assume that $\ell$ has a unique minimum $\tilde{\vtheta}$ on $\T$, then $\argmin_{\T} \widehat{\ell}(\vtheta)=:\widehat{\vtheta} \pto \tilde{\vtheta}$ under some regularity conditions. The regularity conditions might vary across different settings. Consider $\A$ as the parameter to be estimated, the primary aim in the first three steps is to show that $\A$ minimizes the expected loss and establish a uniform convergence result for the random loss function of $\A$.

\subsection{Outline of the second half of the proof}

\textbf{Step 4}: Define $N_{a,b}^j = \sum_{i=1}^N \1\{z_{i,j}^{0}=a\} \1\{\widehat{z}_{i,j}=b\},\ a,b\in[L_j]$ to represent the samples with the wrong local latent class assignments. Derive some upper bounds for the quantities based on $N_{a,b}^j$ by using $\overline{\ell}(\widehat{\A}) - \overline{\ell}(\A^0)$ with the help of the identification assumptions.

\vspace{5pt}
\noindent \textbf{Step 5}: Bound the $\sum_{i=1}^N \1\{ \widehat{\malpha}_i \neq \malpha^0\}$ by using the quantities based on $N_{a,b}^j$ with the help of the discrete structure of the Q-matrix, then obtain the desired classification error rate.

\vspace{5pt}
Assumption~\ref{Ass1} and Assumption~\ref{Ass2} are the regularity conditions for achieving clustering consistency based on the uniform convergence results established in the first half of the proof. We will make more comments about the assumptions in the later proofs.

\subsection{First Half of the Proof of Theorem 1}

\hypertarget{step1}{}
\noindent \textbf{Step 1}. The idea of decomposing $\overline{\ell}(\A) - \ell(\A)$ is to consider \begin{align*}
    \ell(R_{i,j},\widehat{\theta}_{j,z_i}) - \E[\ell(R_{i,j},\overline{\theta}_{j,z_i})] = \left(\ell(R_{i,j},\widehat{\theta}_{j,z_i}) - \ell(R_{i,j}, \overline{\theta}_{j,z_i})\right) + \left( \ell(R_{i,j}, \overline{\theta}_{j,z_i}) - \E[\ell(R_{i,j}, \overline{\theta}_{j,z_i})]\right).
\end{align*}

\noindent The variability in the first term of the right-hand side mainly from the fluctuation in $\abs{\hta - \bta}$, while the randomness in the second term is due to the stochastic nature of $R_{i,j}$.

\hypertarget{Lem1}{}
\begin{Lem}
    Let $(R_{i,j}; 1\leq i\leq N,1\leq j\leq J)$ denote independent Bernoulli trials with parameters $(P_{i,j}; 1\leq i\leq N, 1\leq j\leq J)$. Under a general latent class model, given an arbitrary latent attribute profiles $\A$, there is \begin{align*}
     &\inf_{\T} \E[ \ell(\A,\T|\R)] -\inf_{\T} \ell(\A,\T|\R) \\
     = & \sum_{j=1}^J \sum_{a=1}^{L_j} n_{j,a} (\hta - \bta)^2 + \sum_{i=1}^N \sum_{j=1}^J (P_{i,j} - R_{i,j})(1-2\overline{\theta}_{j,z_i})  \\ = & \sum_{j=1}^J \sum_{a=1}^{L_j} n_{j,a} (\hta - \bta)^2 + \E[X] - X ,\tag{A.7}
    \end{align*}
    \noindent where $X = \sum_{j=1}^J\sum_{i=1}^N R_{i,j}(1- 2\overline{\theta}_{j,z_i})$ is a random variable depending on $\A$ and $L_j$ denotes the number of the distinct local latent classes induced by $\bm{q}_j$ for item $j$.
\end{Lem}

\vspace{5pt}
\noindent \textbf{Proof}. Note $\overline{\ell}(\A) = \sum_i\sum_j(P_{i,j} - \overline{\theta}_{j,z_i})^2 + \sum_i\sum_j P_{i,j}(1 - P_{i,j})$, then \begin{align*}
    &\overline{\ell}(\A) - \ell(\A) \\=& \sum_i\sum_j\left((P_{i,j} - \overline{\theta}_{j,z_i})^2- (R_{i,j} - \widehat{\theta}_{j,z_i})^2  \right) +\sum_i \sum_j P_{i,j}(1 - P_{i,j}) \\ = &  \sum_i\sum_j\left((R_{i,j} - \overline{\theta}_{j,z_i})^2- (R_{i,j} - \widehat{\theta}_{j,z_i})^2  \right)  \\ &+ \sum_i\sum_j\left((P_{i,j} - \overline{\theta}_{j,z_i})^2- (R_{i,j} - \overline{\theta}_{j,z_i})^2  \right) + \sum_{i=1}^N \sum_{j=1}^J P_{i,j}(1-P_{i,j}) \\ = & \sum_i\sum_j\left((R_{i,j} - \overline{\theta}_{j,z_i})^2- (R_{i,j} - \widehat{\theta}_{j,z_i})^2  \right) + \sum_i \sum_j (P_{i,j} - R_{i,j})(1 - 2\overline{\theta}_{j,z_i}).
\end{align*}

\noindent The last equality holds since $R_{i,j}^2 = R_{i,j}$. Given a fixed $\A$, (\hyperlink{A.6}{A.6}) implies that $\sum_{z_{i,j}=a} R_{i,j} = n_{j,a} \hta,\ \sum_{z_{i,j}=a} P_{i,j} = n_{j,a} \bta$, then \begin{align*}
    &\sum_i \sum_j \left( (R_{i,j} - \overline{\theta}_{j,z_i})^2- (R_{i,j} - \widehat{\theta}_{j,z_i})^2\right) \\  = & \sum_{j=1}^J \sum_{a=1}^{L_j}\sum_{z_i =a }  \left( (R_{i,j} - \overline{\theta}_{j,a})^2- (R_{i,j} - \widehat{\theta}_{j,a})^2\right) \\ =& \sum_{j=1}^J \sum_{a=1}^{L_j}\sum_{z_i =a } \left( 2R_{i,j} \hta - 2R_{i,j} \bta + \bta^2 - \hta^2 \right) \\ =& \sum_{j=1}^J \sum_{a=1}^{L_j} \left( 2 n_{j,a} \hta^2 - 2 n_{j,a} \hta \bta + n_{j,a}\bta^2 - n_{j,a} \hta^2  \right) \\ =& \sum_{j=1}^J \sum_{a=1}^{L_j} n_{j,a}(\hta - \bta)^2 .\tag{A.8}
\end{align*}

\noindent This completes the proof of the lemma. \hfill $\square$

\hypertarget{step2}{}
\noindent \textbf{Step 2}. In this step we bound $\sum_j \sum_a n_{j,a} (\hta-\bta)^2$ and $\abs{X-\E[X]}$ separately, for bounding the first term, we have the following lemma. 

\hypertarget{Lem2}{}
\begin{Lem}
    The following event happens with probability at least $1-\delta$, \begin{align*}
        \max_{\A}\left\{ \sum_{j=1}^J \sum_{a=1}^{L_j} n_{j,a} (\hta - \bta)^2\right\} < \frac{1}{2}\left(N \log 2^K + J 2^{K} \log \left( \frac{N}{2^{K}} +1\right) -\log\delta \right).
    \end{align*}
\end{Lem}

\noindent \textbf{Proof}. Under any realization of $\A$, each $\widehat{\theta}_{j,a}$ is an average of $n_{j,a}$ independent Bernoulli random variables $r_{1,j},\dots,r_{N,j}$ with mean $\overline{\theta}_{j,a}$. By applying the Hoeffding inequality, we have \begin{align*}
    \P(\widehat{\theta}_{j,a} \geq \overline{\theta}_{j,a} + t)\leq \exp(-2n_{j,a} t^2),\ \ \P(\widehat{\theta}_{j,a} \leq \overline{\theta}_{j,a} - t)\leq \exp(-2n_{j,a} t^2 ). \tag{A.9}
\end{align*}

\noindent Note that given a fixed $\A$, each $\widehat{\theta}_{j,a}$ can take values only in the finite set $\{0,1/n_{j,a},2/n_{j,a},\dots,1\}$ of cardinality $n_{j,a}+1$. We denote this range of $\widehat{\theta}_{j,a}$ by $\widehat{\Theta}^{j,a}$ and denote the range of the matrix $\widehat{\T} = (\theta_{j,a})$ by $\widehat{\Theta}$. Then $\P(\widehat{\theta}_{j,a} = v) \leq \exp(-2n_a (v- \overline{\theta}_{j,a})^2)$ for any $v\in \widehat{\Theta}^{j,a}$. Since for each of the $J\times 2^K$ entries in $\widehat{\T}$, $\widehat{\theta}_{j,a}$ can independently take on $n_{j,a}+1$ different values, there is $\abs{\widehat{\Theta}} = \prod_j \prod_{a=1}^{L_j} (n_{j,a} + 1)$ with constraint $\sum_{a=1}^{L_j} n_{j,a} = N$. Since $L_j = 2^{K_j}\leq 2^{K}$, we have $\prod_{a=1}^{L_j} (n_{j,a} +1) \leq (1 + N/2^{K})^{2^{K}}$. Denote $\widehat{\Theta}_\epsilon = \{\widehat{\T}\in \widehat{\Theta}: \sum_j\sum_a n_{j,a}(\tilde{\theta}_{j,a} -\overline{\theta}_{j,a})^2\geq \epsilon\}$, then $\widehat{\Theta}_\epsilon \subseteq \widehat{\Theta}$, and 
\begin{align*}
    &\P\left( \sum_j \sum_{a=1}^{L_j} n_{j,a} (\widehat{\theta}_{j,a}-\overline{\theta}_{j,a})^2\geq \epsilon\right) = \sum_{\tilde{\T} \in \widehat{\Theta}_\epsilon} \P\left( \widehat{\T} = \tilde{\T} \right) \\  \leq & \sum_{\tilde{\T} \in \widehat{\Theta}_\epsilon} \prod_j \prod_a \exp\left(-2n_{j,a} (\tilde{\theta}_{j,a}- \overline{\theta}_{j,a})^2\right) \\ =& \sum_{\tilde{\T} \in \widehat{\Theta}_\epsilon}  \exp\left(-2n_{j,a} \sum_j\sum_a (\tilde{\theta}_{j,a}- \overline{\theta}_{j,a})^2\right) \\ \leq & \sum_{\tilde{\T} \in \widehat{\Theta}_\epsilon} \exp(-2\epsilon) \leq \abs{\widehat{\Theta}} e^{-2\epsilon} \\ \leq & \left( \frac{N}{2^{K}} +1 \right)^{J2^{K}} e^{-2\epsilon}. \tag{A.10}
\end{align*}

\noindent The above result holds for any fixed $\A$, we apply a union bound over all the $(2^K)^N$ possible assignments of $\A$ and obtain \begin{align*}\P\left( \max_{\mathbf A} \left\{   \sum_j\sum_a n_{j,a} (\widehat{\theta}_{j,a} - \overline{\theta}_{j,a})^2 \right\} \geq \epsilon\right) \leq 2^{KN}\left(\frac{N}{2^{K}}+1 \right)^{J 2^{K}} e^{-2\epsilon}. \tag{A.11}\end{align*}

\noindent Take $\delta = 2^{KN}\left(\frac{N}{2^{K}}+1 \right)^{J2^{K}} e^{-2\epsilon}$, then $2\epsilon = N\log 2^K + J2^{K} \log(1+N/2^{K}) - \log \delta$. This concludes the proof of Lemma \hyperlink{Lem2}{2}. \hfill $\square$

\hypertarget{Lem3}{}
\begin{Lem}
    Define the random variables $X_{i,j} = R_{i,j} 
 (1-2\overline{\theta}_{j,z_i})$, and $X =\sum_i \sum_j R_{i,j} 
 (1-2\overline{\theta}_{j,z_i})$. Note $X_{i,j}\in[-1,1]$, we apply the Hoeffding's inequality to bound $\abs{X-\E[X]}$ for any realization of $\A$: \begin{align*}
    \P(\abs{X-\E[X]}\geq \epsilon) &\leq 2\exp\left\{ -\frac{2\epsilon^2}{\sum_i\sum_j (1-(-1))^2}    \right\} \\ &\leq 2\exp\left\{ -\frac{\epsilon^2}{2NJ} \right\}. \tag{A.12}
\end{align*} 
\end{Lem}

\noindent With the help of Lemma \hyperlink{Lem2}{2} and Lemma \hyperlink{Lem3}{3}, we next prove the following proposition.

\begin{Prop}\label{Prop1}
    Under the following scaling for some small positive constant $c>0$, \begin{align*}
        \sqrt{J} = O( N^{1- c}),
    \end{align*}
    \noindent we have $\max_{\A} \abs{\overline{\ell}(\A) - \ell(\A)} = o_p(\delta_{N,J})$ where $\delta_{N,J} = N \sqrt{ J} (\log J)^{\tilde{\epsilon}}$ for a small positive $\tilde{\epsilon}>0$. 
\end{Prop}

\noindent \textbf{Proof}. Combining the results of Lemma \hyperlink{Lem2}{2} and Lemma \hyperlink{Lem3}{3}, since there are $(2^K)^N$ possible assignments of $\A$, we apply the union bound to obtain \begin{align*}
    &\P(\max_{\A} \abs{\overline{\ell}(\A) - \ell(\A)} \geq 2\epsilon \delta_{N,J}) \tag{A.13} \\ \leq & (2^K)^N \P\left[  \left\{ \sum_j\sum_a n_{j,a} (\hta-\bta)^2 \geq \epsilon \delta_{N,J}\right\} \cup \{ \abs{X-\E[X]} \geq \epsilon \delta_{N,J}  \}    \right] \\ \leq & \exp\left( N \log(2^K) + J 2^{K} \log\left( \frac{N}{2^{K}} + 1\right)-2\epsilon \delta_{N,J}\right) \\ &+ 2\exp\left( N\log(2^K) - \frac{\epsilon^2 \delta_{N,J}^2}{2NJ}\right).
\end{align*}

\noindent For the second term on the right-hand side of the aforementioned display to converge to zero, we set $\delta_{N,J} = N \sqrt{J} (\log J)^{\tilde{\epsilon}}$ for a small positive constant $\tilde{\epsilon}$. Moreover, under this $\delta_{N,J}$, for the first term to converge to zero as $N,J$ increase, the scaling $\sqrt{J} = O( N^{1- c})$ given in the theorem results in $\P(\max_{\A} \abs{\overline{\ell}(\A) - \ell(\A)} \geq \epsilon \delta_{N,J}) = o(1)$, which implies the result in Proposition \ref{Prop1}.  \hfill $\square$

\hypertarget{step3}{}
\vspace{5pt}
\noindent \textbf{Step 3}. (\hyperlink{A.6}{A.6}) implies that $\overline{\theta}_{j,z_i}^{(\A^0)} = P_{i,j}$, which means that if we plug in the true latent class membership $\A^0$, the estimators will be the corresponding true parameters. According to this property, the following lemma shows that $\A^0$ minimizes the expected loss.

\begin{Lem}\label{Lem4}
    Note $\overline{\ell}(\A^0) - \sum_i\sum_j P_{i,j}(1-P_{i,j}) = \sum_i \sum_j (P_{i,j} - \overline{\theta}_{j, z_i^0})^2 = \sum_i\sum_j(P_{i,j} - P_{i,j})^2=0$, we can obtain \begin{align*}
        \overline{\ell}(\A) - \overline{\ell}(\A^0) = \sum_i \sum_j (P_{i,j} - \overline{\theta}_{j,z_i})^2 \geq 0 .\tag{A.14}
    \end{align*}
\end{Lem}

\noindent Note Lemma \ref{Lem4} holds for any $\A$, it also holds for the estimator $\widehat{\A}$, then \hypertarget{A.15}{}\begin{align*}
    0 \leq \overline{\ell}(\widehat{\A}) - \overline{\ell}(\A^0) = [\overline{\ell}(\widehat{\A}) - \ell(\widehat{\A})] + [\ell(\widehat{\A}) - \ell(\A^0)] + [\ell(\A^0) - \overline{\ell}(\A^0)] . \tag{A.15}
\end{align*}

\noindent Since $\widehat{\A} = \argmin_{\A} \ell(\A)$, we have $\ell(\widehat{\A}) - \ell(\A^0) \leq 0$. Substituting this into (\hyperlink{A.15}{A.15}), we can derive that \begin{align*}
    0\leq \overline{\ell}(\widehat{\A}) - \overline{\ell}(\A^0) \leq 2\sup_{\A} \abs{\overline{\ell}(\A) - \ell(\A)} = o_p(\delta_{N,J}).
\end{align*}

\subsection{Second Half of the Proof of Theorem 1}

\textbf{Step 4}. Motivated by Assumption~\ref{Ass2}, we define $\mathcal{J}:= \{j\in[J];\  \exists k \in[K] \ \mathrm{s.t.}\ \bm{q}_j^0 = \bm{e}_k\}$, which represents the set of all items $j$ that depend on only one latent attribute. Note that $\forall j\in \J$, $\abs{\{\malpha \circ \bm{q}_j^0;\malpha \in \{0,1\}^K\}} = 2$, as $\bm{q}_j$ only contains one required latent attribute, then $\xi(\bm{q}_j^0,\malpha) \in \{1,2\}$ for all $j\in\J$. Without loss of generality, we assume that if $\malpha \circ \bm{q}_j^0 \neq \bm{0}$, then let $\xi(\bm{q}_j^0,\malpha) = 2$, otherwise, let $\xi(\bm{q}_j^0,\malpha) = 1$. We also assume that $\theta_{j,2}^0 > \theta_{j,1}^0, \forall j \in \J$, which aligns with the concept that subjects possessing the required latent attribute tend to perform better. For any $j\in \J$, define \hypertarget{A.16}{}\begin{align*}
    N_{a,b}^j := \sum_{i=1}^N \1\{ z^0_{i,j} = a\} \1\{ \widehat{z}_{i,j} = b\}, \quad (a,b)\in\{1,2\}^2. \tag{A.16}
\end{align*}

\noindent Note $P_{i,j} = \1\{z_{i,j}^0=2\}\theta^0_{j,2} + \1\{z_{i,j}^0=1\}\theta^0_{j,1}$ and $N_{2,2}^j + N_{1,2}^j = \sum_{i=1}^N \1\{ \widehat{z}_{i,j} = 2\}$, $N_{2,1}^j + N_{1,1}^j = \sum_{i=1}^N \1\{ \widehat{z}_{i,j} = 1\}$. By using (\hyperlink{A.6}{A.6}), there is \begin{align*}
    \bttwo  &= \frac{\sum_{i=1}^N \1\{ \widehat{z}_{i,j}=2\} P_{i,j}}{\sum_{i=1}^N \1\{ \widehat{z}_{i,j} = 2\}} \\ &= \frac{\sum_{i=1}^N \1\{ \widehat{z}_{i,j}=2\} (\1\{z_{i,j}^0=2\}\theta^0_{j,2} + \1\{z_{i,j}^0=1\}\theta^0_{j,1})}{\sum_{i=1}^N \1\{ \widehat{z}_{i,j} = 2\}} \\ & = \frac{N_{2,2}^j \theta^0_{j,2} + N_{1,2}^j \theta^0_{j,1}}{N_{2,2}^j + N_{1,2}^j};  \\ \btone & = \frac{N_{2,1}^j \theta^0_{j,2} + N_{1,1}^j \theta^0_{j,1}}{N_{2,1}^j + N_{1,1}^j} .\tag{A.17}
\end{align*}

\noindent Under $\widehat{\A}$, we impose a natural constraint $\bttwo > \btone,\ \forall j\in\J$ on $\widehat{\A}$ for identifiability purpose. This constraint does not change the previous results since $\theta_{j,2}^0 > \theta_{j,1}^0$ allows $\ell(\widehat{\A}) - \ell(\A^0)\leq 0$ in (\hyperlink{A.15}{A.15}) still holds, thus $\overline{\ell}(\widehat{\A}) - \overline{\ell}(\A^0) = o_p(\delta_{N,J})$ still holds under this constraint. Combining $\bttwo > \btone$ and $\theta_{j,2}^0 > \theta_{j,1}^0$, there is \hypertarget{A.18}{}\begin{align*}
    \bttwo > \btone &\Longleftrightarrow (N_{2,2}^j N_{1,1}^j - N_{1,2}^j N_{2,1}^j) \ttwo > (N_{2,2}^j N_{1,1}^j - N_{1,2}^j N_{2,1}^j) \tone \\ & \Longleftrightarrow N_{2,2}^j N_{1,1}^j > N_{2,1}^j N_{1,2}^j .\tag{A.18}
\end{align*}

\noindent From (\hyperlink{A.16}{A.16}), we can obtain \begin{align*}
    \abs{\tone - \btone} = \frac{N^j_{2,1}(\ttwo - \tone)}{N^j_{2,1} + N^j_{1,1}},\quad \abs{\ttwo - \bttwo} = \frac{N^j_{1,2}(\ttwo - \tone)}{N^j_{2,2} + N^j_{1,2}}, \\ \abs{\ttwo - \btone} = \frac{N^j_{1,1}(\ttwo - \tone)}{N^j_{2,1} + N^j_{1,1}}, \quad \abs{\tone - \bttwo} = \frac{N^j_{2,2}(\ttwo - \tone)}{N^j_{2,2} + N^j_{1,2}}.
\end{align*}

\noindent Therefore, \hypertarget{A.19}{}\begin{align*}
    &\ \ \overline{\ell}(\widehat{\A}) - \overline{\ell}(\A^0)\\  = &\  \sum_{j=1}^J \sum_{i=1}^N (P_{i,j} - \overline{\theta}_{j,\widehat{z}_i})^2 \geq \sum_{j\in\J} \sum_{i=1}^N (P_{i,j} - \overline{\theta}_{j,\widehat{z}_i})^2 \\  = &\ \sum_{j\in \J} \left( N_{1,1}^j (\tone - \btone)^2 + N^j_{2,1}(\ttwo  - \btone)^2 +N_{1,2}^j (\tone - \bttwo)^2 + N_{2,2}^j (\ttwo - \bttwo)^2 \right) \\  = &\ \sum_{j\in\J} \left( \frac{N_{1,1}^j (N_{2,1})^2 + N_{2,1}^j(N_{1,1})^2}{(N_{2,1}^j + N_{1,1}^j)^2} + \frac{N_{1,2}^j (N_{2,2})^2 + N_{2,2}^j(N_{1,2})^2}{(N_{2,2}^j + N_{1,2}^j)^2}\right) (\ttwo - \tone)^2 \\  = &\ \sum_{j\in\J} \left( \frac{N_{2,1}^j N_{1,1}^j}{ N_{2,1}^j + N_{1,1}^j} + \frac{N_{2,2}^j N_{1,2}^j}{ N_{2,2}^j + N_{1,2}^j}  \right)(\ttwo - \tone)^2 \\  \geq &\ \delta \sum_{j\in\J} \left( \frac{N_{2,1}^j N_{1,1}^j}{ N_{2,1}^j + N_{1,1}^j} + \frac{N_{2,2}^j N_{1,2}^j}{ N_{2,2}^j + N_{1,2}^j}  \right) \\ \geq &\ \frac{1}{2}\delta\sum_{j\in\J} \left(\min\{N_{2,1}^j,N_{1,1}^j\} + \min \{ N_{2,2}^j , N_{1,2}^j\}\right). \tag{A.19}
\end{align*}

\noindent The second inequality holds since by Assumption~\ref{Ass1}, $(\theta^0_{j,2} - \theta^0_{j,1})^2\geq \delta$. One ideal scenario is that for most $j\in \J$, $\min\{N_{2,1}^j,N_{1,1}^j\} + \min \{ N_{2,2}^j , N_{1,2}^j\} = N_{2,1}^j + N_{1,2}^j = \sum_{i=1}^N \1\{z_{i,j}^0 \neq \widehat{z}_{i,j} \}$, thus the misclassification error for the local latent classes could be bounded relatively tight. The following result confirms this intuition to be accurate.

\begin{Lem}\label{Lem5}
    Define the following random set depending on the estimated latent attribute profiles $\widehat{\A}$ under constraint $\overline{\theta}_{j,2}^{(\A)} > \overline{\theta}_{j,1}^{(\A)},\forall j\in\J$: \begin{align*}
        \J_0 = \{j\in \J; N_{2,1}^j < N_{1,1}^j, N_{1,2}^j < N_{2,2}^j\}; \\
        \J_1 = \{j\in \J; N_{2,1}^j < N_{1,1}^j, N_{1,2}^j > N_{2,2}^j\}; \\
        \J_2 = \{j\in \J; N_{2,1}^j > N_{1,1}^j, N_{1,2}^j < N_{2,2}^j\},
    \end{align*}
    \noindent then under Assumption~\ref{Ass1} and Assumption~\ref{Ass2}, there are $\abs{\J_1} = o_p(\delta_{N,J}/N ),\ \abs{\J_2} = o_p(\delta_{N,J}/N)$
\end{Lem}

\noindent \textbf{Proof}. If $j\in\J_1$, then $\min\{N_{2,1}^j,N_{1,1}^j\} + \min\{N_{2,2}^j, N_{1,2}^j\} = N_{2,1}^j + N_{2,2}^j = \sum_{i=1}^N \1\{z^0_{i,j} = 2 \}$. Under Assumption~\ref{Ass2}, there is $$\sum_{i=1}^N \1\{ z_{i,j}^0 = 2\} \geq N\epsilon.$$

\noindent Then \begin{align*}
    &\  \P\left( \abs{\J_1} \geq b \frac{\delta_{N,J}}{ N \delta}\right) \\ \leq &\ \P\left( \sum_{j\in \J_1} N_{2,1}^j + N_{2,2}^j \geq b \frac{\delta_{N,J}}{N\delta} \cdot N\epsilon \right) \\ \leq &\  \P\left( \overline{\ell}(\widehat{\A}) - \overline{\ell}(\A^0) \geq b\epsilon \delta_{N,J}\right).
\end{align*}

\noindent By noting $b\epsilon$ is a constant and $\overline{\ell}(\widehat{\A}) - \overline{\ell}(\A^0) = o_p(\delta_{N,J})$, then $\abs{\J_1} = o_p(\delta_{N,J}/N )$. Similar arguments give $\abs{\J_2} = o_p(\delta_{N,J}/N)$, which concludes the proof of Lemma \ref{Lem5}. \hfill $\square$

\vspace{8pt}
Note (\hyperlink{A.18}{A.18}) implies that $\min\{N_{2,1}^j,N_{1,1}^j\} + \min \{ N_{2,2}^j , N_{1,2}^j\} \neq N_{1,1}^j +N_{2,2}^j,\ \forall j \in \J$, thus $\J = \J_0 \cup \J_1 \cup \J_2$. The lemma \ref{Lem5} implies that when $\delta_J$ goes to $0$ with a mild rate, the number of elements in $\J_0$ dominates the number of elements in $\J_1 \cup \J_2$, thus for most $j\in\J$, $\min\{N_{2,1}^j,N_{1,1}^j\} + \min \{ N_{2,2}^j , N_{1,2}^j\}$ should be $N_{2,1}^j + N_{1,2}^j = \sum_{i=1}^N \1\{ z_{i,j}^0 \neq \widehat{z}_{i,j} \}$, which represents the number of subjects with the incorrectly assigned local latent classes.

\vspace{8pt}
\noindent \textbf{Step 5}. (\hyperlink{A.19}{A.19}) implies that $o_p(\delta_{N,J}) \geq \overline{\ell}(\widehat{\A}) - \overline{\ell}(\A^0) \geq \delta \sum_{j\in\J_0} \sum_{i=1}^N \1\{ z_{i,j}^0 \neq \widehat{z}_{i,j} \}/2$, next we focus on obtaining a lower bound of $\sum_{j\in \J_0}\sum_{i=1}^N \1\{ z_{i,j}^0 \neq \widehat{z}_{i,j} \}$ to control the classification error rate $N^{-1}\sum_{i=1}^N \1\{\malpha_i^0 \neq \widehat{\malpha}_i \}$.

\vspace{3pt}
Motivated by Assumption~\ref{Ass2}, for each latent attribute $k$, denote $j_k^1$ the smallest integer $j$ such that item $j$ has a $\bm{q}$-vector $\bm{e}_k$, and denote by $j_k^2$ the second smallest integer $j$ such that $\bm{q}_j = \bm{e}_k$, etc. For positive integer $m$, denote \hypertarget{A.20}{}\begin{align*}
    \mathcal{B}^m = \{j_1^m,\dots, j_K^m\}. \tag{A.20}
\end{align*}

\noindent For each $k\in\{1,\dots, K\}$, denote \begin{align*}
    J_{\min}  = \min_{1\leq k\leq K} \abs{\{j\in \J_0;\bm{q}_j^0 = \bm{e}_k\}} , \quad \tilde{J}_{\min} = \min_{1\leq k\leq K} \abs{\{j\in \J;\bm{q}_j^0 = \bm{e}_k\}}.\tag{A.21}
\end{align*}

\noindent Then we have that $\mathcal{B}^m \cap \mathcal{B}^l = \emptyset \ \mathrm{for\ any\ } m\neq l$, thus \hypertarget{A.22}{}\begin{align*}
    & \sum_{i=1}^N \sum_{j\in\J_0} \1\{\xi(\bm{q}_j^0,\malpha_i^0) \neq \xi(\bm{q}_j^0, \widehat{\malpha}_i)\} \\ 
    \geq & \sum_{i=1}^N \sum_{m=1}^{J_{\min}} \sum_{j\in \mathcal{B}^m} \1\{\xi(\bm{q}_j^0,\malpha_i^0) \neq \xi(\bm{q}_j^0, \widehat{\malpha}_i)\} \\ 
    = & J_{\min} \sum_{i=1}^N \sum_{k=1}^K \1\{\xi(\bm{e}_k,\malpha_i^0) \neq \xi(\bm{e}_k, \widehat{\malpha}_i)\} \\
    \geq & J_{\min} \sum_{i=1}^N \1\{ \malpha_i^0 \neq \widehat{\malpha}_i  \}. \tag{A.22}
\end{align*}

\noindent The last inequality holds since $\sum_{k=1}^K \1\{\xi(\bm{e}_k,\malpha_i^0) \neq \xi(\bm{e}_k, \widehat{\malpha}_i)\} \geq \1\{ \malpha_i^0 \neq \widehat{\malpha}_i \}$. Note (\hyperlink{A.22}{A.22}) implies $o_p(\delta_{N,J}/N) \geq J_{\min} N^{-1} \sum_{i=1}^N \1\{ \malpha_i^0 \neq \widehat{\malpha}_i \}$. For simplicity, denote $$\gamma_J = \frac{\delta_{N,J}}{NJ} = \frac{(\log J)^{\tilde{\epsilon}}}{ \sqrt{J}}. $$ 

\noindent Note (\ref{eq12}) in Assumption~\ref{Ass2} implies that $\abs{\J}/J \geq \tilde{J}_{\min}/J \geq \delta_J$ and $J_{\min} \geq \tilde{J}_{\min} - \abs{\J_1 \cup \J_2}$, plug these results into $o_p(\delta_{N,J}/N) \geq J_{\min} N^{-1} \sum_{i=1}^N \1\{ \malpha_i^0 \neq \widehat{\malpha}_i \}$, we can obtain $$o_p\left(\frac{\delta_{N,J}}{N}\right) + \abs{\J_1\cup \J_2} \geq \frac{\tilde{J}_{\min}}{N} \sum_{i=1}^N \1\{\malpha_i^0 \neq \widehat{\malpha}_i\} \geq \frac{J\delta_J}{N} \sum_{i=1}^N \1\{\malpha_i^0 \neq \widehat{\malpha}_i\} .$$

\noindent From Lemma \ref{Lem5} we have $\abs{\J_i} = o_p(\delta_{N,J}/N)$ for $i = 1,2$, which implies that $\abs{\J_1 \cup \J_2} = o_p(\delta_{N,J}/N)$. Plug this into the above inequality, we can conclude that \begin{align*}
    o_p\left( \frac{\delta_{N,J}}{N}\right) \geq \frac{J\delta_J}{N} \sum_{i=1}^N \1\{\malpha_i^0 \neq \widehat{\malpha}_i\} ,
\end{align*}

\noindent which is equivalent to $N^{-1} \sum_{i=1}^N \1\{ \malpha_i^0 \neq \widehat{\malpha}_i  \} = o_p(\gamma_J /\delta_J)$. The proof of the theorem is now complete. \hfill $\square$

The inequality (\ref{eq12}) in Assumption~\ref{Ass2} builds a bridge between the misclassification error for the local latent classes and the misclassification error for the latent attribute profiles $\widehat{\A}$ by using the inequality $\sum_{k=1}^K \1\{\xi(\bm{e}_k,\malpha_i^0) \neq \xi(\bm{e}_k, \widehat{\malpha}_i)\} \geq \1\{ \malpha_i^0 \neq \widehat{\malpha}_i \}$.

\hypertarget{proof_thm_2}{}
\subsection{Proof of Theorem~\ref{Thm2}}

\noindent For notational simplicty, denote $n_{j,a}^0 = \sum_{i=1}^N \1\{z_{i,j}^0 = a \}$. Thus Assumption~\ref{Ass2} implies that $\forall \malpha \in \{0,1\}^K,\ \sum_{i=1}^N \1\{\malpha_i^0 = \malpha \}\geq \epsilon N$ and \begin{align}
    n_{j,a}^0  \geq \frac{2^K}{2^{K_j}} N\epsilon \geq N\epsilon.\tag{A.22}
\end{align}

\noindent Recall that $$\hta = \frac{\sum_{i=1}^N \1\{\widehat{z}_{i,j} = a \} R_{i,j}}{\sum_{i=1}^N \1\{\widehat{z}_{i,j} = a \}}.$$
\noindent Rewrite $\theta_{j,a}^0$ as similar form $$\theta_{j,a}^0 = \frac{\sum_{i=1}^N \1\{{z}_{i,j}^0 = a \} \theta_{j,a}^0}{\sum_{i=1}^N \1\{{z}_{i,j}^0 = a \}} = \frac{\sum_{i=1}^N \1\{{z}_{i,j}^0 = a \} P_{i,j}}{\sum_{i=1}^N \1\{{z}_{i,j}^0 = a \}}.$$

\noindent By triangle inequality we have \begin{align*}
    &\max_{j,a} \abs{\hta - \theta_{j,a}^0} \\ = & \max_{j,a} \abs{\frac{\sum_{i=1}^N \1\{\widehat{z}_{i,j} = a \} R_{i,j}}{\sum_{i=1}^N \1\{\widehat{z}_{i,j} = a \}} - \frac{\sum_{i=1}^N \1\{{z}_{i,j}^0 = a \} P_{i,j}}{\sum_{i=1}^N \1\{{z}_{i,j}^0 = a \}}} \\ \leq & \max_{j,a}\abs{\frac{\sum_{i=1}^N \1\{\widehat{z}_{i,j} = a \} R_{i,j}}{\sum_{i=1}^N \1\{\widehat{z}_{i,j} = a \}} - \frac{\sum_{i=1}^N \1\{\widehat{z}_{i,j} = a \} R_{i,j}}{\sum_{i=1}^N \1\{{z}_{i,j}^0 = a \}}} \\ 
     & + \max_{j,a}\abs{\frac{\sum_{i=1}^N \1\{\widehat{z}_{i,j} = a \} R_{i,j}}{\sum_{i=1}^N \1\{{z}_{i,j}^0 = a \}} - \frac{\sum_{i=1}^N \1\{{z}_{i,j}^0 = a \} R_{i,j}}{\sum_{i=1}^N \1\{{z}_{i,j}^0 = a \}}} \\ & +
    \max_{j,a} \abs{\frac{\sum_{i=1}^N \1\{{z}_{i,j}^0 = a \} R_{i,j}}{\sum_{i=1}^N \1\{{z}_{i,j}^0 = a \}} - \frac{\sum_{i=1}^N \1\{{z}_{i,j}^0 = a \} P_{i,j}}{\sum_{i=1}^N \1\{{z}_{i,j}^0 = a \}}}\\ 
    \equiv & \ \mathcal{I}_1 + \mathcal{I}_2 + \mathcal{I}_3.
\end{align*}

\noindent We then analyze these three terms separately. For the first term, \begin{align*}
    \mathcal{I}_1 &\leq \max_{j,a} \left(\sum_i \1\{\widehat{z}_{i,j} =a \}R_{i,j}\right)\cdot \frac{\sum_i \abs{\1\{\widehat{z}_{i,j}=a\} - \1\{z_{i,j}^0 = a\}}}{n_{j,a}^0 \sum_{i} \1\{ \widehat{z}_{i,j}=a\}}  \\&\leq \max_{j,a}\frac{\sum_i \abs{\1\{\widehat{z}_{i,j}=a\} - \1\{z_{i,j}^0 = a\}}}{n_{j,a}^0} \\ &\leq \frac{1}{\epsilon N}\sum_{i}\1\{\malpha_i^0 \neq \widehat{\malpha}_i\} = o_p\left( \frac{\gamma_J}{\delta_J}\right).
\end{align*}

\noindent The last inequality holds since $\forall j\in[J],\ j\in[L_j],\ \sum_i \abs{\1\{\widehat{z}_{i,j}=a\} - \1\{z_{i,j}^0 = a\}} \leq  \sum_{i}\1\{\malpha_i^0 \neq \widehat{\malpha}_i\}$. For the second term we have \begin{align*}
    \mathcal{I}_2 = \max_{j,a} \frac{\sum_i \abs{R_{i,j} (\1\{\widehat{z}_{i,j}=a\} - \1\{ z_{i,j}^0=a\})}}{n_{j,a}^0} \leq \max_{j,a} \frac{\sum_i \abs{\1\{\widehat{z}_{i,j}=a\} - \1\{z_{i,j}^0 = a\}}}{n_{j,a}^0}.
\end{align*}

\noindent Due to the same reason as $\mathcal{I}_1\pto 0$, we can also conclude that $\mathcal{I}_2 = o_p(\gamma_J/\delta_J)$, thus $\mathcal{I}_1 + \mathcal{I}_2 = o_p(\gamma_J/\delta_J)$. For the third term, we apply Hoeffding's inequality for bounded random variables and obtain $$\P\left( \ \frac{\sum_i \1\{z_{i,j}^0=a\}(R_{i,j} - P_{i,j})}{n_{j,a}^0} \geq t\right) \leq 2\exp(-2 n_{j,a}^0 t^2) \leq 2\exp(-2\epsilon N t^2).$$

\noindent Note the number of $(j,a)$ pair less than or equals to $J\cdot 2^{K}$ under Assumption~\ref{Ass2}, we have for $\forall t>0$, \hypertarget{A.23}{}\begin{align*}
    \P(\mathcal{I}_3\geq t) \leq J 2^{K+1}\exp(-2\epsilon N t^2) .\tag{A.23}
    \end{align*}

\noindent Notably, $2^{K+1}$ remains a constant since $K$ is fixed. By choosing $t = 1/\sqrt{N^{1-\tilde{c}}}$ for a small $\tilde{c}>0$, the tail probability in (\hyperlink{A.23}{A.23}) converges to zero when the scaling condition $\sqrt{J} = O( N^{1- c})$ holds. This implies that $\mathcal{I}_3 = o_p(1/\sqrt{N^{1-\tilde{c}}})$. Bringing together the preceding results, we can conclude that $$\max_{j,a} \abs{\widehat{\theta}_{j,a} - \theta_{j,a}^0} = o_p\left( \frac{\gamma_J}{\delta_J} \right) + o_p\left( \frac{1}{\sqrt{N^{1-\tilde{c}}}}\right).$$

\section{Proof of Theorem~\ref{Thm3}}\label{App2}
The proof shares a similar methodology to that outlined in Appendix \hyperlink{Appendix1}{1} for Theorem~\ref{Thm1}. We will focus on discussing the main differences. The principal distinction between the original GNPC method and the modified GNPC method addressed in Theorem~\ref{Thm1} lies in the imposition of certain parameters to be strictly $0$ or $1$. Specifically, using the ``local latent class'' notation, the constraint (\ref{eq10}) can be reformulated as follows \hypertarget{B.1}{}
\begin{align*}
    &\left\{\xi(\bm{q}_j,\malpha) = 1 \Longleftrightarrow \malpha \cdot \bm{q}_j=0,\ \xi(\bm{q}_j,\malpha) = L_j \Longleftrightarrow \malpha\cdot \bm{q}_j = K_j\right\} ;\\
    &\left\{\theta_{j,1} = 0,\ \theta_{j,L_j}=1\right\}. \tag{B.1}
\end{align*}

\noindent Then $(\widehat{\A},\widehat{\T}) = \argmin_{(\A,\T)} \sum_i \sum_j (R_{i,j} - \theta_{j,z_i})^2$ are defined under the constraint (\hyperlink{B.1}{B.1}) in the following context. The definitions of $\ell(\A,\T|\R)$ and $\overline{\ell}(\A,\T)$ are the same as in (\hyperlink{A.2}{A.2}) and (\hyperlink{A.3}{A.3}). Given any realization of $\A$, denote $\widehat{\T}^{(\A)} = \argmin_{\T} \ell(\A,\T|\R)$ and $\overline{\T}^{(\A)} = \argmin_{\T} \overline{\ell} (\A,\T)$ under the constraint (\hyperlink{B.1}{B.1}). Clearly, (\hyperlink{A.6}{A.6}) still holds for $a\in\{2,\dots,L_j-1\}$ and $\widehat{\theta}_{j,1}^{(\A)} = \overline{\theta}_{j,1}^{(\A)}=0, \ \widehat{\theta}_{j,L_j}^{(\A)} = \overline{\theta}_{j,L_j}^{(\A)}=1$, thus we still have $\E[\widehat{\theta}_{j,a}] = \overline{\theta}_{j,a}$ for all $j\in[J],\ a\in[L_j]$

\vspace{8pt}
\noindent \textbf{Step 1}. This step mirrors Step \hyperlink{step1}{1} in the proof of Theorem~\ref{Thm1}, with the exception of the decomposition form. The item $\sum_{a=1}^{L_j} n_{j,a}(\widehat{\theta}_{j,a} - \overline{\theta}_{j,a})^2$ in Lemma \hyperlink{Lem1}{1} replaced by $\sum_{a=2}^{L_j-1} n_{j,a}(\widehat{\theta}_{j,a} - \overline{\theta}_{j,a})^2$.

\begin{Lem}\label{Lem6}
    Let $(R_{i,j}; 1\leq i\leq N,1\leq j\leq J)$ denote independent Bernoulli trials with parameters $(P_{i,j}; 1\leq i\leq N, 1\leq j\leq J)$. Under a general latent class model, given an arbitrary $\A$, there is \begin{align*}
     &\inf_{\T} \E[ \ell(\A,\T|\R)] -\inf_{\T} \ell(\A,\T|\R) \\
     = & \sum_{j=1}^J \sum_{a=2}^{L_j-1} n_{j,a} (\hta - \bta)^2 + \sum_{i=1}^N \sum_{j=1}^J (P_{i,j} - R_{i,j})(1-2\overline{\theta}_{j,z_i})  \\ = & \sum_{j=1}^J \sum_{a=2}^{L_j-1} n_{j,a} (\hta - \bta)^2 + \E[X] - X, \tag{B.2}
    \end{align*}

    \noindent where $X = \sum_{j=1}^J\sum_{i=1}^N R_{i,j}(1- 2\overline{\theta}_{j,z_i})$ is a random variable depending on $\A$ and $L_j$ denotes the number of local distinct latent classes induced by $\bm{q}_j$ for item $j$.
\end{Lem}

\vspace{5pt}
\noindent \textbf{Proof}. Note $\overline{\ell}(\A) = \sum_i\sum_j(P_{i,j} - \overline{\theta}_{j,z_i})^2 + \sum_i\sum_j P_{i,j}(1 - P_{i,j})$, then we have \begin{align*}
    &\overline{\ell}(\A) - \ell(\A) \\=& \sum_i\sum_j\left((P_{i,j} - \overline{\theta}_{j,z_i})^2- (R_{i,j} - \widehat{\theta}_{j,z_i})^2  \right) +\sum_i \sum_j P_{i,j}(1 - P_{i,j}) \\ = &  \sum_i\sum_j\left((R_{i,j} - \overline{\theta}_{j,z_i})^2- (R_{i,j} - \widehat{\theta}_{j,z_i})^2  \right)  \\ & + \sum_i\sum_j\left((P_{i,j} - \overline{\theta}_{j,z_i})^2- (R_{i,j} - \overline{\theta}_{j,z_i})^2  \right) + \sum_{i=1}^N \sum_{j=1}^J P_{i,j}(1-P_{i,j}) \\ = & \sum_i\sum_j\left((R_{i,j} - \overline{\theta}_{j,z_i})^2- (R_{i,j} - \widehat{\theta}_{j,z_i})^2  \right) + \sum_i \sum_j (P_{i,j} - R_{i,j})(1 - 2\overline{\theta}_{j,z_i}).
\end{align*}

\noindent The last equality holds because $R_{i,j}^2 = R_{i,j}$. Given a fixed $\A$, (\hyperlink{A.6}{A.6}) implies that $\sum_{z_{i,j}=a} R_{i,j} = n_{j,a} \hta,\ \sum_{z_{i,j}=a} P_{i,j} = n_{j,a} \bta$ for $a\in\{2,\dots,L_j-1\}$, then \begin{align*}
    &\sum_i \sum_j \left( (R_{i,j} - \overline{\theta}_{j,z_i})^2- (R_{i,j} - \widehat{\theta}_{j,z_i})^2\right) \\  = & \sum_{j=1}^J \sum_{a=1}^{L_j}\sum_{z_i =a }  \left( (R_{i,j} - \overline{\theta}_{j,a})^2- (R_{i,j} - \widehat{\theta}_{j,a})^2\right) \\ =& \sum_{j=1}^J \sum_{a=2}^{L_j-1}\sum_{z_i =a } \left( 2R_{i,j} \hta - 2R_{i,j} \bta + \bta^2 - \hta^2 \right) \\ =& \sum_{j=1}^J \sum_{a=2}^{L_j-1} \left( 2 n_{j,a} \hta^2 - 2 n_{j,a} \hta \bta + n_{j,a}\bta^2 - n_{j,a} \hta^2  \right) \\ =& \sum_{j=1}^J \sum_{a=2}^{L_j-1} n_{j,a}(\hta - \bta)^2 .\tag{B.3}
\end{align*}

\noindent The second equality holds since for $a\in\{1,L_j\},\ \widehat{\theta}_{j,a} = \overline{\theta}_{j,a}$ thus $(R_{i,j} - \overline{\theta}_{j,z_i})^2- (R_{i,j} - \widehat{\theta}_{j,z_i})^2=0$ when $z_i \in\{1,L_j\}$. This concludes the proof of Lemma \ref{Lem6}. \hfill $\square$

\vspace{8pt}
\noindent \textbf{Step 2}. This step closely parallels Step \hyperlink{step2}{2} in the proof of Theorem~\ref{Thm1}. Given that $\sum_{a=2}^{L_j-1}n_{j,a}(\widehat{\theta}_{j,a} - \overline{\theta}_{j,a})^2\leq \sum_{a=1}^{L_j}n_{j,a}(\widehat{\theta}_{j,a} - \overline{\theta}_{j,a})^2$, and note that $\overline{\theta}_{j,a}$ still falls within $[0,1]$ under the constraint (\hyperlink{B.1}{B.1}), we can directly apply Lemma \hyperlink{Lem2}{2} and Lemma \hyperlink{Lem3}{3} to control $\overline{\ell}(\A) - \ell(\A)$. The same scaling condition and error rate outlined in Proposition \ref{Prop1} remain valid.

\hypertarget{step3_2}{}
\vspace{8pt}
\noindent \textbf{Step 3}. The main difference between this step and Step \hyperlink{step3}{3} in the proof of Theorem~\ref{Thm1} derives from a fact that $\A^0$ does not necessarily minimize $\overline{\ell}(\A)$. As a result, the statement that $\overline{\ell}(\A)\geq \overline{\ell}(\A^0)$ for any $\A$ may not be true. However, for $a\in\{2,\dots,L_j-1 \}$, the condition $\overline{\theta}_{j,a}^{(\A^0)} = P_{i,j}$ continues to be valid, thus we define \hypertarget{B.4}{}\begin{align*}
    S_{N,J}:&=\overline{\ell}(\A^0) - \sum_i \sum_j P_{i,j}(1-P_{i,j}) = \sum_i \sum_j (P_{i,j} - \overline{\theta}_{j,z_i^0})^2 \\ &= \sum_j \sum_{a=1} \sum_{z_i^0=a} (P_{i,j}-0)^2 + \sum_j \sum_{a=L_j} \sum_{z_i^0=a} (1-P_{i,j})^2 . \tag{B.4}
\end{align*}

\noindent By using the same argument in (\hyperlink{A.15}{A.15}), the following inequality holds 
\begin{align*}
    -S_{N,J} \leq \overline{\ell}(\widehat{\A}) - \overline{\ell}(\A^0) \leq 2\sup_{\A} \abs{\overline{\ell}(\A) - \overline{\ell}(\A)} = o_p(\delta_{N,J}).
\end{align*}

\noindent \textbf{Step 4}. Inspired by Assumption~\ref{Ass2}, we define $\mathcal{J}:= \{j\in[J]; \exists k \in[K] \ \mathrm{s.t.}\ \bm{q}_j^0 = \bm{e}_k\}$, representing the set of all items $j$ that are only dependent on one latent attribute. With the constraint (\hyperlink{B.1}{B.1}) in place, when $j\in\mathcal{J}$, it holds that $\overline{\theta}_{j,a} \in\{0,1\},\forall a\in\{1,2\}$. Correspondingly, we define $\N_{a,b}^j:=\{i\in[N]; z_{i,j}^0 = a, \widehat{z}_{i,j}=b \}$ and $N_{a,b}^j = |\N_{a,b}^j|$ in the same manner as in (\hyperlink{A.16}{A.16}), yielding

\hypertarget{B.5}{}
\begin{align*}
    \overline{\ell}(\widehat{\A}) - \overline{\ell}(\A^0) &= \sum_{j=1}^J\sum_{i=1}^N (P_{i,j} - \overline{\theta}_{j,z_i})^2 - S_{N,J} \\ &\geq \sum_{j\in\mathcal{J}}\sum_{i=1}^N (P_{i,j} - \overline{\theta}_{j,z_i})^2 - S_{N,J} \\ &= \sum_{j\in\mathcal{J}}\left(\sum_{i \in \N_{1,1}^j} P_{i,j}^2 + \sum_{i\in\N_{1,2}^2} (P_{i,j}-1)^2 + \sum_{i\in\N_{2,1}^j} P_{i,j}^2 + \sum_{i\in\N_{2,2}^j}(P_{i,j}-1)^2\right) - S_{N,J} \\ &\geq \sum_{j\in \mathcal{J}}\left( \sum_{i\in\N_{1,2}^2} (P_{i,j}-1)^2 + \sum_{i\in\N_{2,1}^j} P_{i,j}^2\right) - S_{N,J} \\ &> \frac{1}{4} \sum_{j\in\mathcal{J}}(N_{1,2}^j + N_{2,1}^j) - S_{N,J} \\ &= \frac{1}{4}\sum_{j\in\mathcal{J}} \sum_{i=1}^N \1\{ z_{i,j}^0 \neq \widehat{z}_{i,j} \} - S_{N,J}. \tag{B.5}
\end{align*}

\noindent The last inequality holds since by Assumption~\ref{Ass3}, $1-P_{i,j}>1/2$ for $i\in\N_{1,2}^j$ and $P_{i,j}>1/2$ for $i\in\N_{2,1}^j$. By the definition of $\lambda_{N,J}^2$ in (\ref{eq17}), we have $S_{N,J} \leq NJ \lambda_{N,J}^2$. Combining with (\hyperlink{B.5}{B.5}), the following inequality holds: \hypertarget{B.6}{}\begin{align*}
    (4NJ)\lambda_{N,J}^2 + o_p(\delta_{N,J}) \geq \sum_{j\in\mathcal{J}} \sum_{i=1}^N\1\{ z_{i,j}^0 \neq \widehat{z}_{i,j} \}. \tag{B.6}
\end{align*}

Defining $\tilde{J}_{\min} = \min_{1\leq k\leq K} \abs{\{j\in \J;\bm{q}_j^0 = \bm{e}_k\}}$, note that Assumption~\ref{Ass2} implies $\abs{\mathcal{J}}/J \geq \tilde{J}_{\min}/J \geq \delta_J$. By using the same definition of $\mathcal{B}^m$ in (\hyperlink{A.20}{A.20}) and note $\mathcal{B}^m\cap \mathcal{B}^{l} \neq \emptyset$ for any $m\neq l$, then \hypertarget{B.7}{}\begin{align*}
    &\ \sum_{j\in\mathcal{J}} \sum_{i=1}^N \1\{ z_{i,j}^0 \neq \widehat{z}_{i,j} \} \\
     \geq &\ \sum_{i=1}^N \sum_{m=1}^{\tilde{J}_{\min}} \sum_{j\in\mathcal{B}^m} \1\{ \xi(\bmq_j^0,\malpha_i^0) \neq \xi(\bmq_j^0 , \widehat{\malpha}_i)\}\\ 
    = &\ \tilde{J}_{\min}\sum_{i=1}^N \sum_{k=1}^K \1\{\xi(\bm{e}_k,\malpha_i^0) \neq \xi(\bm{e}_k,\widehat{\malpha}_i) \} \\
    \geq &\ \tilde{J}_{\min} \sum_{i=1}^N \1\{ \malpha_i^0 \neq \widehat{\malpha}_i \} \\ \geq &\ (J\delta_J)\cdot  \sum_{i=1}^N \1\{ \malpha_i^0 \neq \widehat{\malpha}_i \}. \tag{B.7}
\end{align*}

\noindent The first inequality holds since $\cup_{m=1}^{\tilde{J}_{\min}} \mathcal{B}^m \subseteq \mathcal{J}$ and the second inequality holds due to $\sum_{k=1}^K \1\{\xi(\bm{e}_k,\malpha_i^0) \neq \xi(\bm{e}_k, \widehat{\malpha}_i)\} \geq \1\{ \malpha_i^0 \neq \widehat{\malpha}_i \}$. Combining (\hyperlink{B.6}{B.6}) with (\hyperlink{B.7}{B.7}) we can conclude that \begin{align*}
    \frac{4\lambda_{N,J}^2}{\delta_J} + o_p\left( \frac{\gamma_J}{\delta_J}\right) \geq \frac{1}{N} \sum_{i=1}^N \1\{ \malpha_i^0 \neq \widehat{\malpha}_i \}.
\end{align*}
This  completes the proof of Theorem~\ref{Thm3}. \hfill $\square$

There is a significant distinction between Theorem~\ref{Thm1} and Theorem~\ref{Thm3}. According to the current literature on the consistency results for the NPC method and the GNPC method, having the true class membership minimize the expected loss is crucial for establishing clustering consistency. In the discussion of the GNPC method, the true latent attribute profiles $\A^0$ might not minimize the expected loss function $\overline{\ell}$ if we permit certain parameters to be zero and one. Denote $\tilde{\A} = \argmin_{\A} \overline{\ell}(\A)$, then $\widehat{\A}$ might closer to $\tilde{\A}$ rather than $\A^0$. 

$\lambda_{N,J}^2$ represents the average squared distance between the true parameters and the parameters constrained to 0 or 1. From Step \hyperlink{step3_2}{3}, we observe that $(NJ)\lambda_{N,J}^2\geq \overline{\ell}(\A^0) - \overline{\ell}(\tilde{\A})\geq 0$, implying that if $\lambda_{N,J}^2$ is small, then $\A^0$ nearly minimizes the expected loss $\overline{\ell}$. From this observation, $\lambda_{N,J}^2$ can be interpreted as the cost incurred by forcing some parameters to be exactly zero and one. The result given in Example \ref{Eg1} somewhat supports this idea.

\hypertarget{proof_eg_1}{}
\section{Proof of Example 1}\label{sup_sec_example1}
Under the conditions given in Example \ref{Eg1}, we can derive that $(P_{i,j} - \overline{\theta}_{j,z_i^0})^2 = (1/2 - \lambda)^2$. This is due to the fact that under the true latent class profiles $\A^0$, the pairing of $(P_{i,j} ,\overline{\theta}_{j,z_i^0})$ can either be $(1/2 - \lambda , 0)$ or $(1/2+\lambda,1)$, thereby validating the equation (\ref{eq20}).

For the plug-in profiles $\A^1$ as constructed in the example, first consider the $2M$ samples assigned to $\bm{e}_1$ in $\A^1$. {Given that $\bm{e}_1$ contains just one latent attribute, as per the Q-matrix, when $j\in\{1,2\}$, $\bm{e}_1$ contains half of the necessary latent attributes for items $j$. Therefore the true memberships for the samples assigned to $\bm{e}_1$ are half $(1,1,1,1)$—with corresponding true parameter $(\lambda+1/2)$—and half $(0,0,0,0)$, with corresponding true parameter $(1/2 -\lambda)$. Thus, by (\hyperlink{A.6}{A.6}), we know 
\begin{equation*}
    \overline{\theta}_{j,\bm{e}_1}^{(\A^1)} = \frac{\sum_{i=1}^N \1\{z_{i,j}^{(\A^1)} = \bm{e}_1\}P_{i,j}}{\sum_{i=1}^N \1\{z_{i,j}^{(\A^1)}=\bm{e}_1\}} = \frac{1}{2}(1/2+\lambda) + \frac{1}{2}(1/2-\lambda) = \frac{1}{2},
\end{equation*}
when $j\in\{1,2\}$. Next, note that $\bm{e}_1$ does not include any of the required latent attributes for items $3$ and $4$. Therefore $\overline{\theta}_{j,\bm{e}_1}^{(\A^1)}=0$ for $j \in \{3,4\}$.}

Denote $\malpha_i$ as the row vectors of $\A^1$ and $\malpha_i^0$ as the row vectors of $\A^0$, then \hypertarget{B.8}{}
\begin{align*}
    &\sum_{\malpha_i = \bm{e}_1}\sum_{j=1}^4 \left(P_{i,j} - \overline{\theta}_{j,\bm{e}_1}^{(\A^1)}\right)^2 \\  
    = & \sum_{\malpha_i = \bm{e}_1} \left( \left(P_{i,1} - \frac{1}{2} \right)^2 + \left(P_{i,2} - \frac{1}{2} \right)^2+ P_{i,3}^2  +P_{i,4}^2 \right) \\ =&
    \sum_{\malpha_i = \bm{e}_1, \malpha_i^0 = (1,1,1,1)} \left( \left(P_{i,1} - \frac{1}{2} \right)^2 + \left(P_{i,2} - \frac{1}{2} \right)^2+ P_{i,3}^2  +P_{i,4}^2 \right) \\ & +  \sum_{\malpha_i = \bm{e}_1, \malpha_i^0 = (0,0,0,0)} \left( \left(P_{i,1} - \frac{1}{2} \right)^2 + \left(P_{i,2} - \frac{1}{2} \right)^2+ P_{i,3}^2  +P_{i,4}^2 \right) \\ = & \sum_{\malpha_i =\bm{e}_1, \malpha_i^0 = (1,1,1,1)} \left( 2\lambda^2+ 2\left( \frac{1}{2}+\lambda\right)^2 \right) + \sum_{\malpha_i =\bm{e}_1, \malpha_i^0 = (0,0,0,0)} \left( 2\lambda^2+ 2\left( \frac{1}{2}-\lambda\right)^2 \right) \\ = & M\cdot \left( 2\lambda^2+ 2\left( \frac{1}{2}+\lambda\right)^2 \right) + M\cdot \left( 2\lambda^2+ 2\left( \frac{1}{2}-\lambda\right)^2 \right)\\ = &8M\cdot\left(\lambda^2 + \frac{1}{8} \right)\nonumber\\ =& \frac{1}{4}(NJ) \left( \lambda^2 + \frac{1}{8}\right) .\tag{B.8}
\end{align*}
{Here the first equality is obtained by plugging in the expressions of $\overline{\theta}_{j,\bm{e}_1}^{(\A^1)}$, and the last equality holds since $M = N/8$ and $J = 4$.}

Since the aforementioned analysis for $\bm{e}_1$ is applicable to the other $\bm{e}_i,i\in\{2,3,4\}$ as well, we can insert the derived result in (\hyperlink{B.8}{B.8}) into the equation (\hyperlink{B.4}{B.4}). Consequently, we can deduce that the equation (\ref{eq21}) indeed holds true. \hfill $\square$

\section{Additional Simulation Results}\label{supp_sec_simu}
\begin{landscape}
\begin{figure}[htbp!]
\centering    
\subfigure{
        \includegraphics[width=9.2in,height=1.2in]{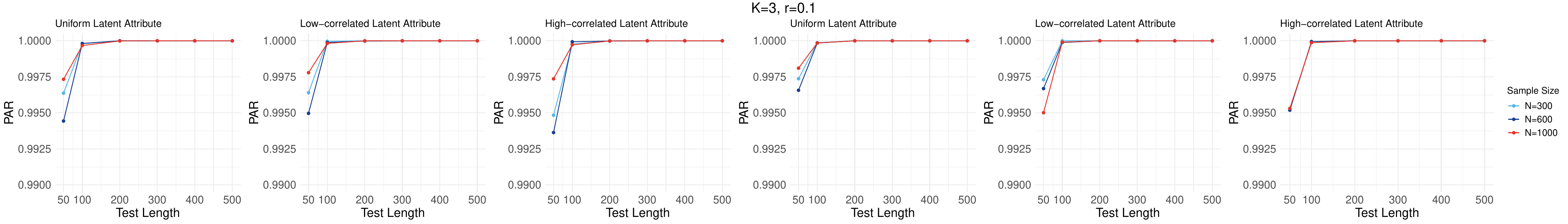}}
        \subfigure{
        \includegraphics[width=9.2in,height=1.2in]{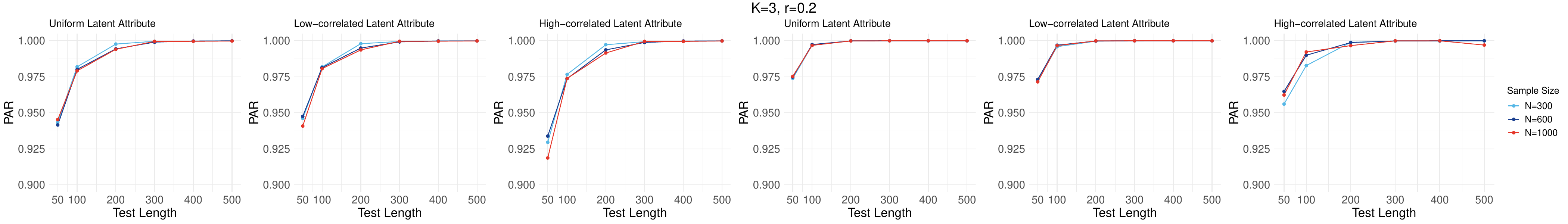}}
        \subfigure{
        \includegraphics[width=9.2in,height=1.2in]{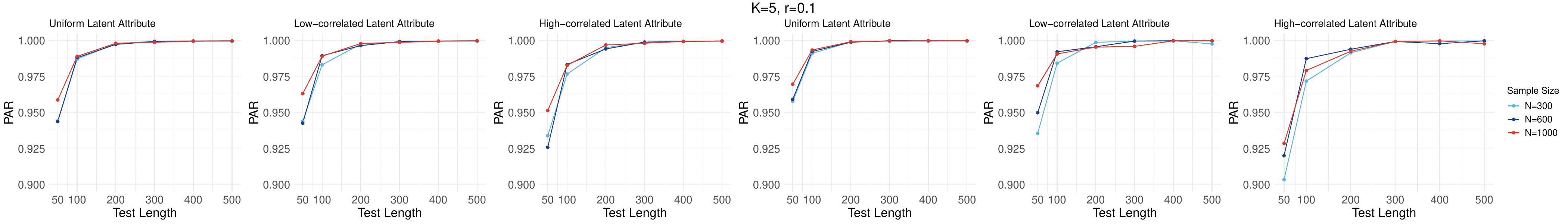}}
        \subfigure{
        \includegraphics[width=9.2in,height=1.2in]{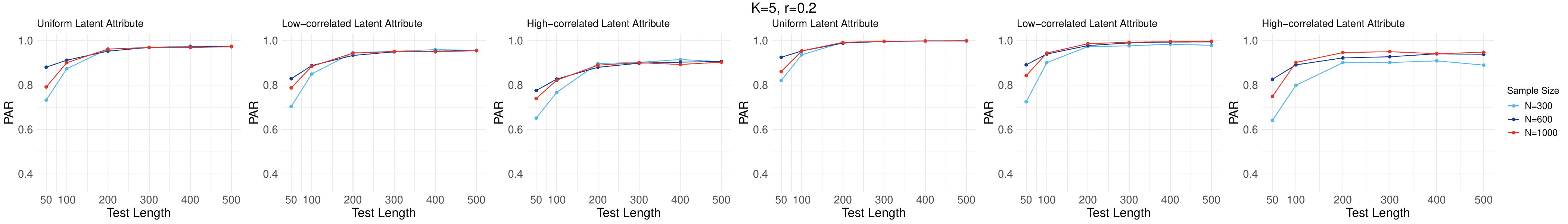}}
    \caption{{PARs for data generated under the DINA model with random initialization. In each row, the left three subfigures show the results of the original GNPC method, while the right three subfigures show the results of the modified GNPC method. }}
    \label{random:DINA}
\end{figure}
\end{landscape}

\begin{landscape}
\begin{figure}[htbp!]
\centering    
\subfigure{
        \includegraphics[width=9.2in,height=1.2in]{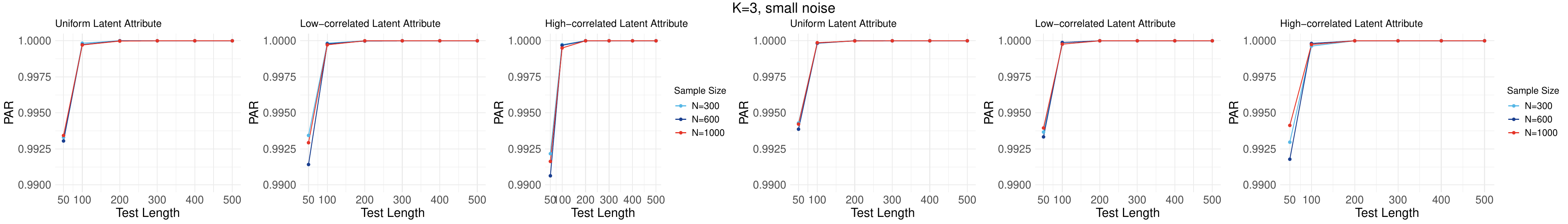}}
        \subfigure{
        \includegraphics[width=9.2in,height=1.2in]{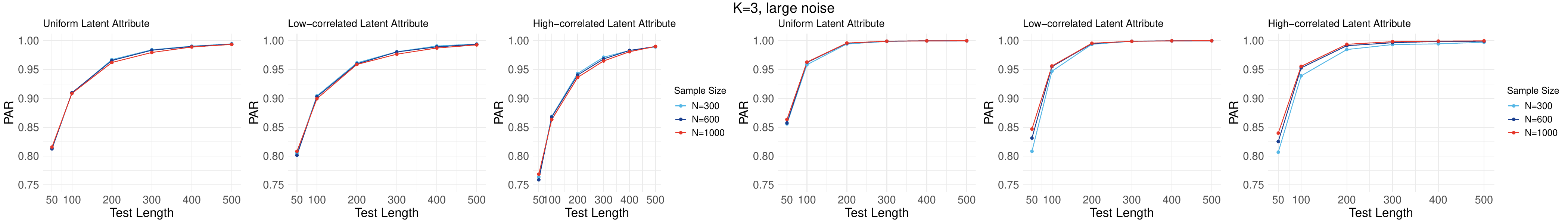}}
        \subfigure{
        \includegraphics[width=9.2in,height=1.2in]{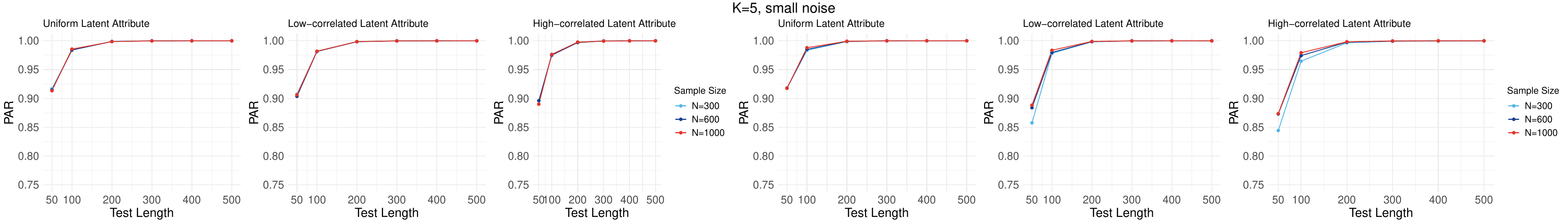}}
        \subfigure{
        \includegraphics[width=9.2in,height=1.2in]{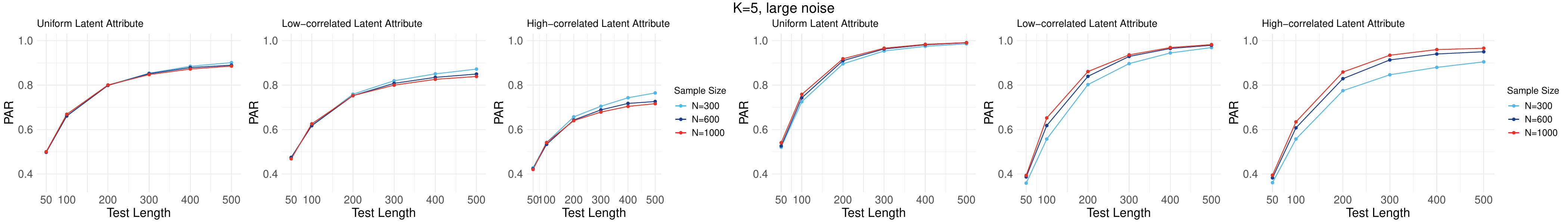}}
    \caption{{PARs for data generated under the GDINA model with random initialization. In each row, the left three subfigures show the results of the original GNPC method, while the right three subfigures show the results of the modified GNPC method. }}
    \label{random:GDINA}
\end{figure}
\end{landscape}

In the main text, the simulation studies present estimation results using both the original GNPC method~\cite{chiu2019consistency} and the modified GNPC method~\cite{Ma2023}, where the latent attributes estimated by the NPC method are used as initial values. Empirically, the NPC method provides a good initial estimate for the proficiency classes, which accelerates the convergence of both GNPC methods. Yet our theoretical results demonstrate that consistent initialization may not be necessary. 
In this section, we replicate the simulation studies under all settings from the main text but randomly generate the initial values for the proficiency classes of all subjects. To mitigate the risk of local optima, we repeat the procedure five times and select the estimation with the smallest loss, $\ell(\hat\A, \hat\T|\R)$, as the final estimate. The estimation results are shown in Figure~\ref{random:DINA} and Figure~\ref{random:GDINA} for data generated under the DINA model and the GDINA model, respectively.

In general, both the original GNPC method and the modified GNPC method with random initialization achieve comparable estimation errors across all scenarios to those initialized using the NPC method. These results further demonstrate that, in large-scale measurement, both the original and modified GNPC methods can consistently estimate the latent attributes without relying on NPC initialization.

\begin{landscape}
\begin{figure}[htbp!]
\centering    
\subfigure{
        \includegraphics[width=9.2in,height=1.2in]{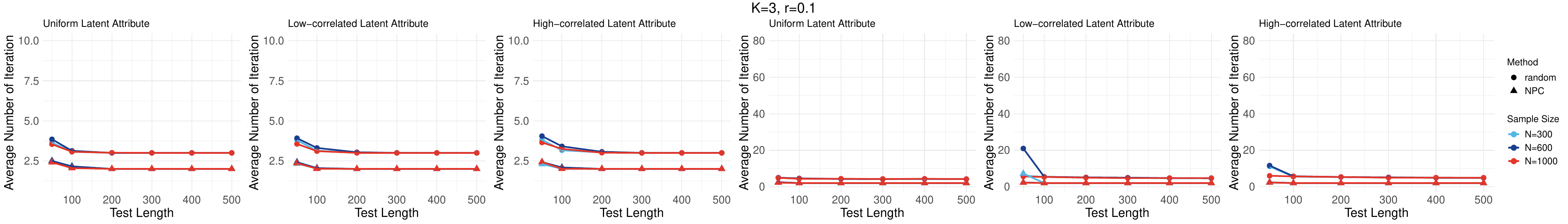}}
        \subfigure{
        \includegraphics[width=9.2in,height=1.2in]{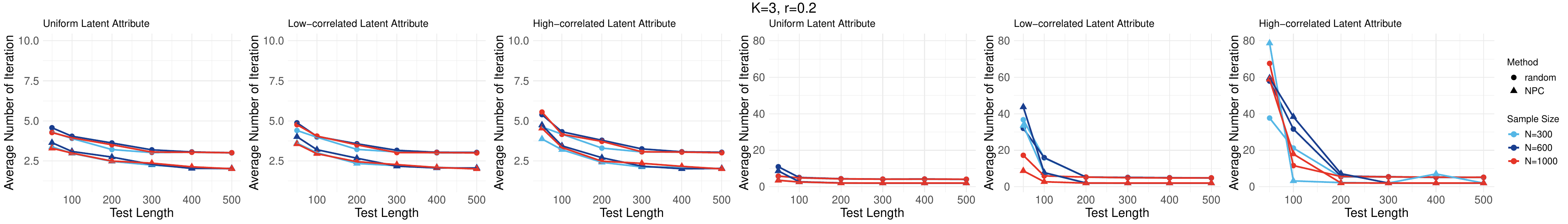}}
        \subfigure{
        \includegraphics[width=9.2in,height=1.2in]{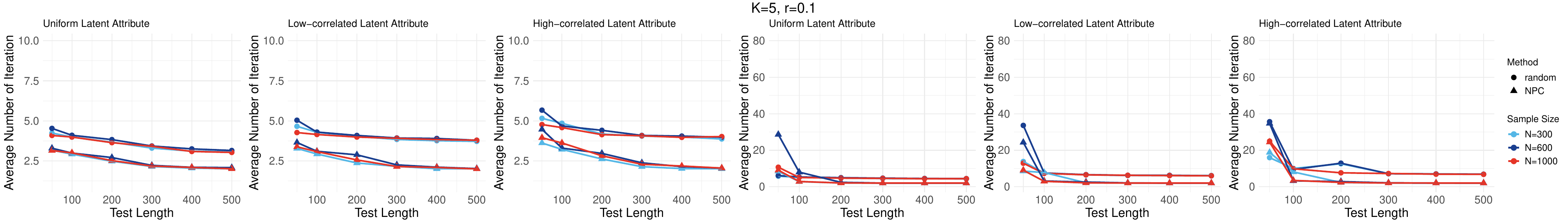}}
        \subfigure{
        \includegraphics[width=9.2in,height=1.2in]{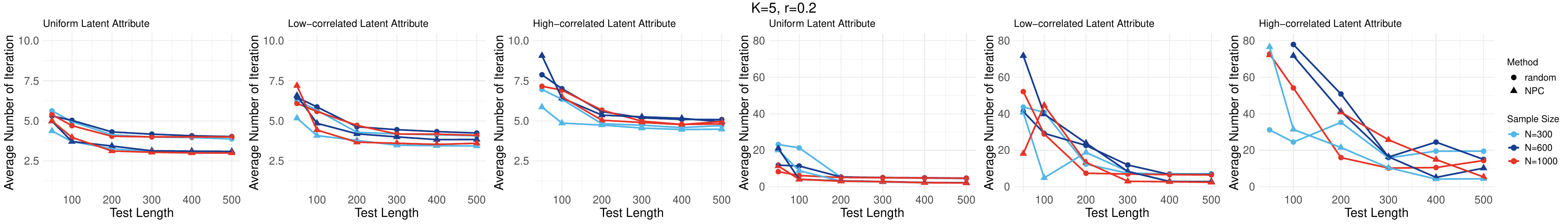}}
    \caption{{Average number of iteration of GNPC using random initialization and NPC initialization. Data are generated under the DINA model. In each row, the left three subfigures show the results of the original GNPC method, while the right three subfigures show the results of the modified GNPC method.}}
    \label{compare:DINA}
\end{figure}
\end{landscape}

\begin{landscape}
\begin{figure}[htbp!]
\centering    
\subfigure{
        \includegraphics[width=9.2in,height=1.2in]{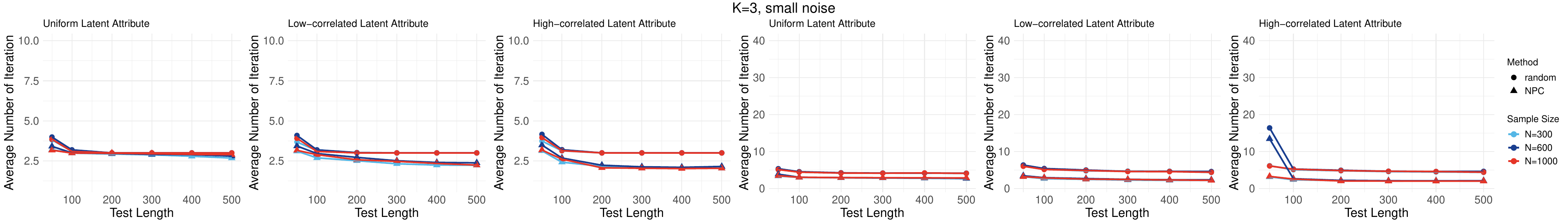}}
        \subfigure{
        \includegraphics[width=9.2in,height=1.2in]{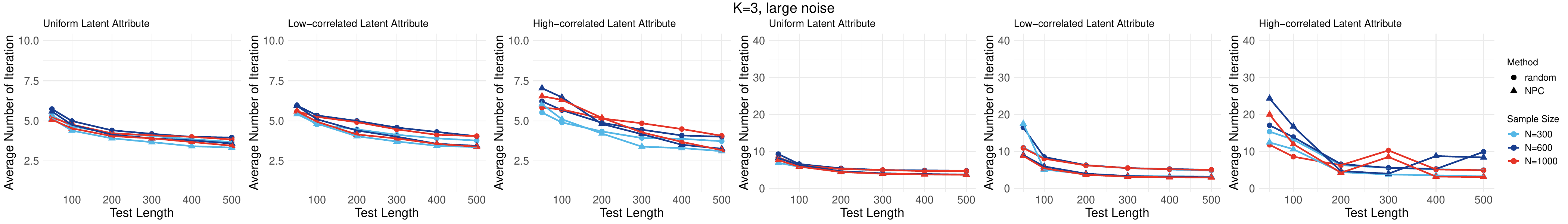}}
        \subfigure{
        \includegraphics[width=9.2in,height=1.2in]{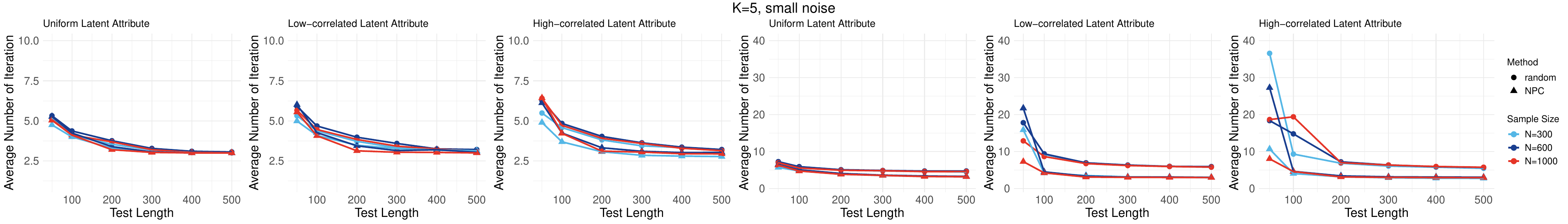}}
        \subfigure{
        \includegraphics[width=9.2in,height=1.2in]{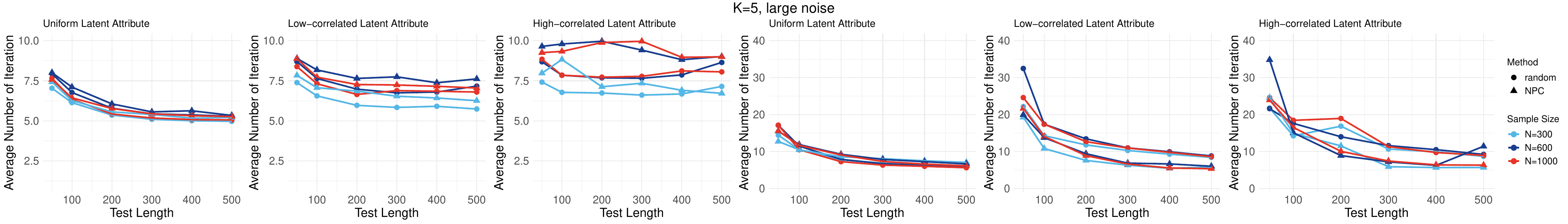}}
    \caption{{Average number of iteration of GNPC using random initialization and NPC initialization. Data are generated under the GDINA model. In each row, the left three subfigures show the results of the original GNPC method, while the right three subfigures show the results of the modified GNPC method.}}
    \label{compare:GDINA}
\end{figure}
\end{landscape}

We also examine the increase in computational cost for both GNPC methods when using random initialization. Specifically, we compare the average number of iterations required by both methods under random initialization versus initialization using the NPC method. The results are shown in Figure~\ref{compare:DINA} and Figure~\ref{compare:GDINA} for data generated under the DINA model and the GDINA model, respectively. Our findings indicate a noticeable increase in computational cost with random initialization. Additionally, in more complex scenarios—such as those with higher noise levels or a greater number of attributes $(K=5)$—the modified GNPC method requires a significantly larger number of iterations to achieve convergence. Overall, these observations highlight that the GNPC methods could consistently estimate latent attributes even with random initialization, but greater computational efficiency can be achieved when initialized with a well-chosen estimate.

\end{document}